\documentclass[leqno,12pt]{article}

\usepackage{amsthm,amssymb}
\usepackage{amsbsy,amsfonts,amsmath}

\makeatletter
      
       \@addtoreset{equation}{section}
\makeatother

\newcommand{\bysame}{\leavevmode\hbox to
3em{\hrulefill}\,}

\newcommand{\al}{{\alpha}}
\newcommand{\calE}{{\cal E}}
\newcommand{\calEtimes}{{\calE^\times}}
\newcommand{\EPiakg}{{{\cal D}(\calE,\Pi,a,k,g)}}
\newcommand{\EPiako}{{{\cal D}(\calE,\Pi,a,k,0)}}
\newcommand{\PiB}{{\Pi^B}}
\newcommand{\Pical}{{\Pi_c(\al)}}
\newcommand{\Rkg}{{R(k,g)}}

\newcommand{\BBbb}{\Bbb}

\newcommand{\BBox}{\Box}

\newcommand{\Der}{{{\rm D}^1}}

\begin{document}

      \begin{center}
{\bf  {A finite number of defining relations
and a UCE theorem of the
elliptic Lie algebras
and superalgebras
with
rank $\geq 2$}}
\end{center}

Hiroyuki Yamane
\par
Department of Pure and Applied Mathematics, Graduate
School of Information
Science
and Technology, Osaka University, Toyonaka 560-0043, Japan
\par
e-mail:yamane@ist.osaka-u.ac.jp
\newline\newline
\begin{abstract}
In this paper, we give a finite number of defining
relations satisfied by a finite number of generators
for the elliptic Lie algebras and superalgebras
${\frak g}_R$ with
rank $\geq 2$. Here the $R$'s denote
the reduced and non-reduced elliptic root systems with
rank $\geq 2$.
We also show that
if ${\cal L}$ is an extended
affine Lie algebra (EALA) whose non-isotropic roots form
the $R$,
then there exists a natural homomorphism ${\cal F}:{\frak
g}_R
\rightarrow{\cal L}$, which also give a
universal central extension (UCE) surjective map from 
$[{\frak g}_R,{\frak g}_R]$ to the core of ${\cal L}$.
(More precisely, we take a ${\bar {\frak g}}_R$
instead of the ${\frak g}_R$.)
\end{abstract}

\section*{Introduction}

In 1985, K.~Saito \cite{saito1} introduced the notion of
the (reduced and non-reduced) {\it extended affine root
systems}.
In this paper, we also call them the
SEARS's.  Let $R$ be an SEARS. Let ${\cal V}:=\BBbb{R}R$
and
${\cal V}^0:=\{v\in{\cal V}|s_\al(v)=v\,\mbox{for all
$\al\in R$}\}$,
where $s_\al$ denotes the reflection with respect to an
$\al$.
Let $m:=\dim{\cal V}^0$
and $l:=\dim{\cal V}-m$. We say that the $m$ is the {\it
nullity} of the
$R$
and say that the $l$ is the {\it rank} of the $R$.
The $R$ is {\it reduced} if $\BBbb{R}\al\cap
R=\{\al,-\al\}$
for all $\al\in R$.
We notice that
if $m=0$, the $R$ is a finite root system, and that if
$m=1$,
the $R$ is an affine root system. If $m=2$, the $R$ is
called
an {\it elliptic root system} (or ERS for short); see also
Subsec.~1.3.
(The notation ${\cal V}$ and ${\cal V}^0$ is used only in
Introduction.)
\par
In 1997, B.~Allison, S.~Azam, S.~Berman, Y.~Gao and
A.~Pianzola
\cite{Allisonetal} introduced and studied root systems
defined by different
axioms
from those of the reduced SEARS's.
They also studied Lie algebras ${\cal L}=({\cal L},
(\,,\,),{\cal H})$
associated with their
root systems. The ${\cal L}$ is called the extended affine
Lie
algebra
(EALA for short) (see \cite{Allisonetal}).
In 2002, S.~Azam
\cite{Azam} showed
that there
exists a natural one-to-one correspondence between their
root systems and
the reduced
SEARS's.
\par
Let $R$ be an ERS. Let ${\widetilde
\pi}:{\cal V}\rightarrow{\cal V}/{\cal V}^0$
be a natural projective map. The $R$ is called {\it
simply-laced}
if $l\geq 2$ and ${\widetilde \pi}(R)$ is
a simply-laced finite root system. In 2000, K.~Saito and
D.~Yoshii
\cite{saitoyoshii} studied a
Lie algebra ${\frak g}_R$ whose non-isotropic roots form a
simply-laced ERS
$R$,
showed that the ${\frak g}_R$ is isomorphic to the
(2-variable)
toroidal Lie algebra \cite{MoodyRaoYokonuma}
of type ${\rm ADE}$,
and gave a Serre-type theorem for the ${\frak g}_R$
(see also the next paragraph). We also notice the results
in \cite{miki1}, \cite{miki2}, \cite{takebayashi} for the
${\frak g}_R$ with the ${\widetilde \pi}(R)$ of type ${\rm
A}$.
\par
In this paper,
for every ERS
$R$ with
$l\geq 2$, we
give a Serre-type theorem for a Lie algebra or
superalgebra ${\frak g}_R$
having the property
that the
system formed by its non-isotropic roots is isomorphic to
the $R$ (see
Theorems~\ref{theorem:main} and
\ref{theorem:ellipticgenerators}); in other
words, we
give a finite number of defining relations of the ${\frak
g}_R$ satisfied
by Chevalley generators.
The
${\frak g}_R$ is not a Lie algebra but a Lie superalgebra
if and only if the $R$ is not reduced.
(In the text, the
${\frak g}_R$ shall be denoted as the
${\frak g}_{{\cal D}}$.)
We call them the {\it elliptic Lie (super)algebra}.
If $R$ is simply-laced, our Serre-type theorem
(Theorem~\ref{theorem:ellipticgenerators}) coincides with
K.~Saito and
D.~Yoshii's Serre-type theorem \cite{saitoyoshii}.
\par
Let $R$ be an ERS. If
there exists a one dimensional subspace $G_{\BBbb{R}}$ of
${\cal V}^0$
such that $G_{\BBbb{R}}\cap\BBbb{Z}R\ne\{0\}$ and the $\pi
(R)$ is a
reduced
affine root system, where
$\pi:{\cal V}\rightarrow{\cal V}/G_{\BBbb{R}}$ is the
natural projection,
then we call the $R$ the {\it reduced marked
ERS}. (More precisely, in the text, we call the pair
$(R,G)$ the reduced marked
ERS (see Subsec.~1.3), where
$G=G_{\BBbb{R}}+\sqrt{-1}G_{\BBbb{R}}$.)
If the $R$ is a reduced marked
ERS with $l\geq 2$, then the Serre-type theorem for the
${\frak g}_R$ has
already
been given by the author \cite{yamane2} in 2004.
\par
In the process of giving the Serre-type theorem,
we also give a Saito-type classification theorem of all
the ERS's $R$
with $l\geq 2$
(see Theorems~\ref{theorem:general} and
\ref{theorem:submain}).
To prove the classification theorem, we use K.~Saito's
classification
theorem
\cite{saito1} of the reduced marked ERS's, and use the
${\frak g}_R$
to show that each $R$ in the classification theorem
really exists.
\par We also give universality theorems of the ${\frak
g}_R$
with $l\geq 2$.
(see Theorems~\ref{theorem:invarianttheorem},
\ref{theorem:qinvarianttheorem} and
\ref{theorem:centralextension}).
Especially, we see that if $R$ is a reduced ERS with
$l\geq 2$ and
if
there exists a basis $\{\delta,\,a\}$ of ${\cal V}^0$
such that $(\al+{\cal V}^0)\cap
R=\al+\BBbb{Z}\delta+\BBbb{Z}a$
for every $\al\in R$,
then the ${\Der}({\frak g}_R):=
[{\frak g}_R,{\frak g}_R]$ is isomorphic to the
universal central
extension (UCE for short) of the Lie algebra ${\frak
f}_R\otimes\BBbb{C}[t_1^{\pm
1},t_2^{\pm
1}]$, where the ${\frak f}_R$ is a finite dimensional
simple Lie algebra
whose root system is isomorphic to the ${\widetilde
\pi}(R)$, i,e., the ${\Der}({\frak g}_R)$ is
a (2-variable) toroidal
Lie algebra in the sense of  \cite{MoodyRaoYokonuma}
(see Corollary~\ref{corollary:MEY}).
We also treat the quantum
tori elliptic Lie algebras studied in
\cite{BermanGaoKrylyuk}
(see Subsec.~5.2).
We show that if an $({\cal L},(\,,\,),{\cal H})$ is an
EALA whose non-isotropic roots form an ERS $R$ with $l\geq
2$,
there exists a natural homomorphism
${\cal F}:{\frak g}_R\rightarrow{\cal L}$
(More precisely,
we take a ${\bar {\frak g}}_R=
{\bar {\frak g}}_{{\cal D}}$
instead of the ${\frak g}_R$).
It turns out that
the ${\cal F}_{|{\Der}({\frak g}_R)}
:{\Der}({\frak g}_R)\rightarrow {\cal F}({\Der}({\frak
g}_R))$ 
is a UCE, and
the image ${\cal F}({\Der}({\frak g}_R))$
is the core of the $({\cal L},(\,,\,),{\cal H})$
(see Corollary~\ref{corollary:MEY})
(see \cite[Chap.~I Definition 2.20]{Allisonetal} for the
term).
\par This paper is organized as follows. Notice that we
always
assume
$l\geq 2$.
In \S1, we discuss the ERS's $R=\Rkg$. In
\S2 we give a definition of the elliptic Lie
algebras and superalgebras ${\frak g}_R={\frak
g}_{{\cal D}}$
and give Theorem~\ref{theorem:main}, which state a
root space decomposition of the ${\frak g}_R$,
and
Theorem~\ref{theorem:submain}, which give a Saito-type
classification of all the ERS's $R$ with $l\geq 2$. In
\S3, we give
a proof of Theorem~\ref{theorem:main} after supposing
Lemma~\ref{lemma:Cartannondeg}. In $\S4$, we give a proof
of Lemma~\ref{lemma:Cartannondeg}.
In \S5, we give Theorem~\ref{theorem:invarianttheorem},
which states that the ${\frak g}_R$
is the maximal ones among the Lie (super)algebras having
the root space $R$ and satisfying some additional
conditions,
and give Theorem~\ref{theorem:centralextension},
which is a UCE theorem of the
${\Der}({\frak g}_R)$.
In \S6, we give Theorem~\ref{theorem:ellipticgenerators},
which
is a natural extension of
K.~Saito and
D.~Yoshii's Serre-type theorem \cite{saitoyoshii}
of the simply-laced elliptic Lie algebras.

\section{Preliminary}
\subsection{Pre-elliptic base system}

In this paper, we set $\BBbb{Z}_+:=\{r\in\BBbb{Z}|r\geq
0\}$, i.e.,
$\BBbb{Z}_+=\{0\}\cup\BBbb{N}$. We also set
$\BBbb{Z}_-:=\{r\in\BBbb{Z}|r\leq 0\}$.
\par
Let $l$ be a fixed positive integer.
Through out this paper, we assume
$l\geq 2$.
Let $\calE$ be an $l+4$-dimensional
$\BBbb{C}$-vector space.
Let $J:\calE\times \calE\rightarrow\BBbb{C}$
be a non-degenerate symmetric bilinear form.
Let $\calEtimes:=\{x\in\calE|J(x,x)\ne 0\}$.
If $x\in\calEtimes$,
we call $x$ a {\it non-isotropic element},
let $x^\vee:={\frac {2x} {J(x,x)}}$
and define $s_x\in{\rm GL}(\calE)$
by $s_x(y)=y-J(x^\vee,y)x$.
Let
$\{\al_0,\ldots,\al_l,\,\Lambda_\delta,\,a,\,\Lambda_a\}$
be a basis of $\calE$ satisfying the following. \par
({\rm B}1) The $(l+1)\times(l+1)$-matrix
$A:=(J(\al_i^\vee,\al_j))_{0\leq i,j\leq l}$
is an affine type generalized Cartan matrix
\cite{kac1}, \cite{kac2}.
Then the $A$ is called
${\bf A}_l^{(1)}$ ($l\geq 2$), ${\bf B}_l^{(1)}$ ($l\geq
3$),
${\bf C}_l^{(1)}$ ($l\geq 2$), ${\bf D}_l^{(1)}$ ($l\geq
4$),
${\bf E}_l^{(1)}$ ($l=6,\,7.\,8$), ${\bf F}_4^{(1)}$
($l=4$),
${\bf G}_2^{(1)}$ ($l=2$),
${\bf A}_{2l}^{(2)}$ ($l\geq 2$), ${\bf A}_{2l-1}^{(2)}$
($l\geq 3$),
${\bf D}_{l+1}^{(2)}$ ($l\geq 2$), ${\bf E}_6^{(2)}$
($l=4$)
or ${\bf D}_4^{(3)}$ ($l= 2$).
(See \cite[Tables 1-4]{kac2}
and \cite[\S4.8 TABLE Aff 1,2,3]{kac1}.)
(The $A$ is neither ${\bf A}_1^{(1)}$
nor  ${\bf A}_2^{(2)}$ since $l\geq 2$.)
The numbering of the elements $\al_0,\ldots,\al_l$
of $\Pi$ is the same as in \cite[\S4.8 TABLE Aff 1, 2,
3]{kac1}.
\par
({\rm B}2) $J(\Lambda_\delta,\Lambda_\delta)=0$,
$J(a,\al_i)=0$,
$J(a,\Lambda_\delta)=0$,
$J(a,a)=0$,
$J(\Lambda_a,\al_i)=0$,
$J(\Lambda_a,a)=1$ and $J(\Lambda_a,\Lambda_a)=0$ and
\begin{equation}\label{eqn:lambdazeroalzero}
J(\Lambda_\delta,\al_i)={\frac {1} {r}}\delta_{i0},
\end{equation} where $r=2$ if $A={\bf A}_{2l}^{(2)}$;
and, otherwise, $r=1$. \par
Let $\Pi:=\{\al_0,\ldots,\al_l\}$.
Let $W$ be the subgroup of ${\rm GL}(\calE)$
generated by $s_\al$ ($\al\in\Pi$), i.e., $W$
is the affine Weyl group.
\par For a subset $S$ of $\Pi$, let $W_S$ be the subgroup
of $W$ generated by $s_\al$ ($\al\in S$).
Let $l(w)$ be the length of $w\in W$
with respect to $s_\al$ ($\al\in\Pi$).

\newtheorem{lemma}{Lemma}[section]
\begin{lemma} \label{lemma:walgamma}
Let $\al$, $\beta\in\Pi$ and $w\in W$ be such that
$w(\al)=\beta$ and $l(w)>0$. Then there exists
a $\gamma\in\Pi\setminus\{\al\}$ and a $w^\prime\in
W_{\{\al,\gamma\}}$
such that $w^\prime(\al)\in\{\al,\gamma\}$ and
$l(w^\prime)+l(w(w^\prime)^{-1})=l(w)$.
\end{lemma}

This can be proved by a well-known argument (see
\cite[Proof
of Proposition 8.20]{Jantzen}).
\newline\par
A function $f:\Pi\rightarrow \BBbb{C}$ is called
{\it $W$-invariant} if $f(\al)=f(\beta)$
for every $(\al,\beta)\in\Pi\times\Pi$
with $\beta=w(\al)$ for some $w\in W$.
By Lemma~\ref{lemma:walgamma}, we see the following.

\begin{lemma} \label{lemma:falfbeta}
Keep the notation as above.
A function $f:\Pi\rightarrow\BBbb{C}$ is $W$-invariant
if and only if $f(\al)=f(\beta)$ for every
$(\al,\beta)\in\Pi\times\Pi$
with $J(\al^\vee,\beta)=J(\al,\beta^\vee)=-1$.
\end{lemma}

Let $k:\Pi\rightarrow \BBbb{N}$ be a $W$-invariant
function such that
${\rm G.C.D.}\{k(\al)|\al\in\Pi\}=1$.
Let $g:\Pi\rightarrow 2^\BBbb{Z}$ be a $W$-invariant
function, where $2^\BBbb{Z}$ is the power set of
$\BBbb{Z}$, i.e., the set of the subsets of $\BBbb{Z}$.
We call a quintuple ${\cal D}=\EPiakg$
of such $\calE$, $\Pi$, $a$, $k$ and $g$ a
{\it {pre-elliptic base system}} (PEBS for short).
If $g(\al)=\emptyset$ for every $\al$, the $g$ is also
denoted
by $0$. \par
For $x\in\calE$ and a subset $B$ of $\BBbb{C}$,
let $Bx:=\{bx\in\calE|b\in B\}$;
moreover, for a subset $X$ of $\calE$, let
$BX:=\sum_{x\in X}Bx$.
(If $B$ is an empty set $\emptyset$, then $BX=\emptyset$.)
For subsets
$S$ and $T$ of $\calE$, let $S+T:=\{x+y\in\calE|x\in
S,\,y\in T\}$;
if $T=\{x\}$, let $x+S:=T+S$.
(If $S=\emptyset$, $S+T=\emptyset$.)
\par
Let ${\cal D}$ be a PEBS.
Let
\begin{equation}
\Rkg:=\bigcup_{w\in W}w\Bigl(
\bigcup_{\al\in\Pi}\bigl(
(\al+\BBbb{Z}k(\al)a)\cup(2\al+g(\al)k(\al)a)
\bigr)\Bigr).
\label{eqn:Rkg}
\end{equation}
Then
\begin{equation}
\Rkg\subset\Bigl((\BBbb{Z}_+\Pi+\BBbb{Z}a)
\cup(\BBbb{Z}_-\Pi+\BBbb{Z}a)\Bigr)\setminus\BBbb{Z}a.
\label{eqn:alzkala}
\end{equation}
\begin{lemma}\label{lemma:zsetimes}
Keep the notation as above.
Let $S$ be a subset of $\Pi$.
Let $\lambda\in (\BBbb{Z}S)\cap{\cal E}^\times$. Then
there exists
a $w\in W_S$ such that
\begin{equation}
w(\lambda)\in \Bigr(\bigcup_{\al\in S}\BBbb{Z}\al\Bigl)
\cup
\Bigr(\BBbb{Z}S\setminus(\BBbb{Z}_+S\cup\BBbb{Z}_-S)\Bigl).
\nonumber
\end{equation}
\end{lemma} \par {\it Proof.}
By the definition of $\Pi$ in (B1), there exists an
$e\in\BBbb{C}\setminus\{0\}$ such that
$eJ(\al_i,\al_j)\in\BBbb{R}$ ($0\leq i,\,j\leq l$),
and $eJ(\al_i,\al_i)>0$ ($0\leq i\leq l$).
Then the symmetric bilinear form
$(eJ)_{|\BBbb{R}\Pi\times\BBbb{R}\Pi}:\BBbb{R}\Pi\times\BBbb{R}\Pi
\rightarrow\BBbb{R}$
is semipositive definite.
For $\mu=\sum_{\al\in\Pi}b_\al\al\in\BBbb{R}\Pi$ with
$b_\al\in\BBbb{R}$,
let $ht(\mu):=\sum_{\al\in\Pi}b_\al\in\BBbb{R}$. \par
Let $\lambda$ be as in the statement. We may assume
$\lambda\in(\BBbb{Z}_+S\cup\BBbb{Z}_-S)\setminus\{0\}$.
Moreover we may assume $\lambda\in\BBbb{Z}_+S$.
We use an induction on $ht(\lambda)\in\BBbb{N}$.
If $ht(\lambda)=1$, then $\lambda\in S$.
We assume $ht(\lambda)>1$.
We may assume $\lambda\notin\BBbb{N}\al$ for any $\al\in
S$.
Since $\lambda\in{\cal E}^\times$,
$eJ(\lambda,\lambda)>0$.
Hence there there exists an $\al\in S$ such that
$eJ(\lambda,\al)>0$.
Notice that $s_\al(\lambda)
=\lambda-{\frac {2eJ(\lambda,\al)} {eJ(\al,\al)}}\al$.
If
$s_\al(\lambda)\notin\BBbb{Z}S\setminus(\BBbb{Z}_+S\cup\BBbb{Z}_-S)$,
then $s_\al(\lambda)\in\BBbb{Z}_+S$ and
$ht(s_\al(\lambda))<ht(\lambda)$.
This completes the proof.
\hfill $\BBox$ \newline\par
For a subset $S$ of $\Pi$, let
\begin{equation}\nonumber
\Rkg_S:=\Rkg\cap((\oplus_{\al\in S}\BBbb{C}\al)\oplus
\BBbb{C}a).
\end{equation}
\begin{lemma}\label{lemma:Rkgs}
Keep the notation as above. Then
\begin{equation}
\Rkg_S=\bigcup_{w\in W_S}w\Bigl(
\bigcup_{\al\in S}\bigl(
(\al+\BBbb{Z}k(\al)a)\cup(2\al+g(\al)k(\al)a)
\bigr)\Bigr).
\nonumber
\end{equation}
In particular, for $\al\in\Pi$, we have
\begin{equation}
\Rkg_{\{\al\}}=
\bigcup_{\epsilon\in \{1,-1\}}\bigl(
(\epsilon\al+\BBbb{Z}k(\al)a)\cup(2\epsilon\al+g(\al)k(\al)a).
\nonumber
\end{equation}

\end{lemma}
{\it Proof.} If $|S|=1$, then the lemma follows
immediately from the
definition.
     From this, together with (\ref{eqn:alzkala}) and
Lemma~\ref{lemma:zsetimes},
the lemma for a general $S$ follows; notice that $w(a)=a$
for all
$w\in W_S$.
\hfill $\BBox$

\subsection{Elliptic and quasi-elliptic base
systems}
Here we introduce the notions of an elliptic base system
(EBS for short)
and a quasi-elliptic base system (QEBS for short).
In Theorem~\ref{theorem:submain}, we shall show that these
notions are
equivalent.
In Theorem~\ref{theorem:general}, we shall show how
the
EBS's
are associated
with the elliptic root
systems. \par
Let ${\cal D}$ be an PEBS with $l\geq 2$.
For $\al\in\Pi$, let
\begin{equation}
\Pical:=\{\beta\in\Pi|\beta\ne\al,\,
J(\beta,\al)\ne 0\}.
\nonumber
\end{equation}
Define a subset $\PiB$ of $\Pi$ by
\begin{equation}
\PiB:=\{\al\in\Pi|\forall \beta\in\Pical,\,\,
J(\al^\vee,\beta)=-2\}.
\nonumber
\end{equation}
We call ${\cal D}$ an {\it an quasi elliptic base system}
(QEBS for short) if the following hold.
\newline\par
({\rm KG}1) If $\al\in\Pi$, $\beta\in\Pical$ and
$J(\beta^\vee,\al)=-1$, then
${\frac {k(\beta)} {k(\al)}}\in\BBbb{Z}$ and
$J(\al^\vee,\beta){\frac {k(\al)} {k(\beta)}}\in\BBbb{Z}$.
\par
({\rm KG}2) $g(\al)=\emptyset$ if $\al\notin\PiB$. \par
({\rm KG}3) If $\al\in\PiB$ and $\beta\in\Pical$, then
$g(\al){\frac {k(\al)} {k(\beta)}}=
\emptyset$, $\BBbb{Z}$, $2\BBbb{Z}$ or $2\BBbb{Z}+1$.
\newline\par
We call a PEBS ${\cal D}$ an {\it elliptic base system}
(EBS for short)
if the following holds.
\begin{equation}
\forall \al\in \Rkg,\quad s_\al(\Rkg)=\Rkg.
\nonumber
\end{equation}

\begin{lemma} \label{lemma:kg}
Let ${\cal D}$ a PEBS with $l\geq 2$. If ${\cal D}$
is an EBS, then it is also a QEBS.

\end{lemma} \par
{\it Proof.}  The axiom ({\rm KG}1) follows from
Lemma~\ref{lemma:Rkgs} and
the
following ({\it
cf.}
\cite[Proof of (6.1) Assertion]{saito1}).
\begin{equation}
\left\{\begin{array}{c}
s_\al
s_{\al+mk(\al)a}(\beta)=\beta+mJ(\al^\vee,\beta)k(\al)a,\\
s_\beta s_{\beta+mk(\beta)a}(\al)=\al+mk(\beta)a
\end{array}\right.
\nonumber
\end{equation}
for $m\in\BBbb{Z}$.
\par
If $\beta\in\Pical$, then $\BBbb{Z}\owns
J((2\al)^\vee,\beta)
={\frac {J(\al^\vee,\beta)} {2}}$. Hence, if
$g(\al)\ne\emptyset$, then
$\al\in\PiB$, which implies the axiom ({\rm KG}2). \par
Let $\al\in\PiB$. Assume $g(\al)\ne\emptyset$.
Let $n\in\BBbb{Z}$ be such that
$2\al+nk(\al)a\in\Rkg$.
Then
\begin{equation}
s_\al s_{\al\pm k(\al)a}(2\al+nk(\al)a)=
2\al+(n\mp 4)k(\al)a
\nonumber
\end{equation} and
\begin{equation}
s_\al s_{2\al+nk(\al)a}(2\al+nk(\al)a)=
2\al-nk(\al)a\,.
\nonumber
\end{equation}
Hence $g(\al)=\BBbb{Z}$, $2\BBbb{Z}$, $2\BBbb{Z}+1$,
$4\BBbb{Z}$ or
$4\BBbb{Z}+2$.
Let $\beta\in\Pical$.
Then
\begin{equation}
s_\beta
s_{\beta+mk(\beta)a}(2\al+nk(\al)a)=2\al+(nk(\al)-2mk(\beta))a
\nonumber
\end{equation} Hence $g(\al)=\BBbb{Z}$, $2\BBbb{Z}$ or
$2\BBbb{Z}+1$
if $k(\beta)=k(\al)$.
Moreover
\begin{equation}
s_\al
s_{2\al+nk(\al)a}(\beta+mk(\beta)a)=\beta+(mk(\beta)-nk(\al))a\,.
\nonumber
\end{equation} Hence $g(\al)=2\BBbb{Z}$, $4\BBbb{Z}$ or
$4\BBbb{Z}+2$
if $k(\beta)=2k(\al)$. This implies the axiom ({\rm KG}3)
and completes the
proof.
\hfill
$\BBox$
\newline\par
Converse of Lemma~\ref{lemma:kg} shall be given in
Theorem~\ref{theorem:submain}.
\par If $\EPiako$ is a QEBS, i.e., $g=0$,
then it is called a {\it {special}} QEBS
(SQEBS for short).

\subsection{Elliptic root systems}

Keep the notation in \S1. Notice that $l\geq 2$. For a
subset
$S$ of ${\cal E}$, let $S^\perp:=\{x\in\calE|\forall y\in
S,\,J(x,y)=0\}$.
Following \cite{saito1} (and \cite{saitoyoshii}), we say
that
a subset $R$
of ${\cal E}$ is an {\it {elliptic root system}}
(ERS for short)
of rank $l$
if it satisfies the following. \newline\par
({\rm SER}1)\quad
$\forall x,\,\forall y\in \BBbb{R}R,\,\,
J(x,x)J(y,y)\in\BBbb{R}_+$, \par
({\rm SER}2)\quad
$\dim_{\BBbb{C}}(\BBbb{C}R\cap R^\perp)=2$, \par
({\rm SER}3)\quad
$\dim_{\BBbb{C}}\BBbb{C}R=l+2={\rm
rank}_{\BBbb{Z}}\BBbb{Z}R$, \par
({\rm SER}4)\quad
$\forall \al\in R,\,\,s_\al(R)=R$,\par
({\rm SER}5)\quad
$\forall \al,\,\forall \beta\in
R,\,\,J(\al^\vee,\beta)\in\BBbb{Z}$, \par
({\rm SER}6)\quad
If $R=R_1\cup R_2$, $R_2\subset (R_1)^\perp$, then
$R_1\ne\emptyset$ or
$R_2\ne\emptyset$,
\newline\newline
where $\BBbb{R}_+=\{x\in\BBbb{R}|x\geq 0\}$. \par
Let $R$ be an ERS. We call the ${\cal E}$ for the $R$ the
{\it
base space}. A one dimensional subspace $G$ of
$\BBbb{C}R\cap R^\perp$
is called a {\it marking line} if $G\cap \BBbb{Z}R\ne
\{0\}$.
The pair $(R,G)$
of the above $R$ and $G$ is called a {\it marked elliptic
root
system}
(MERS for short). An MERS $(R,G)$ is called a {\it reduced
marked elliptic
root
system} (RMERS for short) if
\begin{equation}
\forall \al,\,\forall \beta\in R,\,\,2\al-\beta\notin G\,.
\end{equation}

By \cite{saito1}, we have the following.

\newtheorem{theorem}{Theorem}[section]
\begin{theorem}[(6.4) of \cite{saito1}]
\label{theorem:saitotheorem} If $(R,G)$ be an RMERS,
then there exists an SQEBS $\EPiako$ such that
$R=R(k,0)$, $G=\BBbb{C}a$ and ${\cal E}$ is the base space
of $R$.
If $\EPiako$ is an SQEBS, then $(R(k,0),\BBbb{C}a)$ is an
RMERS;
in particular, $\EPiako$ is an EBS.
\end{theorem}
      In \cite{saito1}, we do not need the assumption
$l\geq
2$
for Theorem~\ref{theorem:saitotheorem}.

\begin{theorem}
\label{theorem:general} (1) If ${\cal D}$ is an EBS, then
$(R(k,g),
\BBbb{C}a)$
is an MERS. \par
(2) Let $(R,G)$ be an MERS.
Then there exists an EBS ${\cal D}$ such that
$R=R(k,g)$, $G=\BBbb{C}a$ and ${\cal E}$ is the base space
of $R$.
\end{theorem}
\par {\it Proof.} The statement (1) is clear.
We prove the statement (2). Let
\begin{equation}
R^\prime:=\{\al\in R\,|\,{\frac \al 2}\notin R +G\}.
\nonumber
\end{equation}
Let $R^{\prime\prime}:=R\setminus R^\prime$. We have
\begin{equation}
\label{eqn:fraconetwoRprimeprimesubsetRprime}
{\frac 1 2}R^{\prime\prime}\subset R^\prime +G
\end{equation}
because, if $\al\in R^{\prime\prime}$ is such that
${\frac \al 2}\notin (R^\prime +G)$,
then ${\frac \al 4}\in R +G$ and
$J(\al^\vee, {\frac \al 4})={\frac 1 2}\notin\BBbb{Z}$,
contradiction.
\par
We show \begin{equation} \label{eqn:QRprimeZequalQRZ}
\BBbb{Z}R^\prime=\BBbb{Z}R.
\end{equation}
Clearly $\BBbb{Z}R^\prime\subset \BBbb{Z}R$ holds.
Let $\beta\in R^{\prime\prime}$. By
(\ref{eqn:fraconetwoRprimeprimesubsetRprime}),
      $\beta=2\al+x$ for some $\al\in R^\prime$
and some $x\in G$. Notice that
\begin{equation}
\forall\gamma\in R,\,\,\sigma_\gamma(R^\prime)=R^\prime
\quad{\rm
and}\quad\sigma_\gamma(R^{\prime\prime})=R^{\prime\prime}\,.
\end{equation}
It follows that $R^\prime\owns$
$\sigma_\al\sigma_\beta(\al)=$
$\sigma_\al(\al-\beta)=$ $-\al-(\beta-4\al)=$
$3\al-\beta=$ $\al-x$.
Hence $x\in \BBbb{Z}R^\prime$. Hence $\BBbb{Z}R\subset
\BBbb{Z}R^\prime$,
as desired. \par
By (\ref{eqn:QRprimeZequalQRZ}), we see that
the $(R^\prime,G)$
is an RMERS whose base space is ${\cal E}$.
By Theorem~\ref{theorem:saitotheorem}, there exists
an SQEBS
$\EPiako$ such that $R(k,0)=R^\prime$ and $G=\BBbb{C}a$.
It follows from Lemma~\ref{lemma:Rkgs} that
\begin{equation}\label{eqn:alGcapR}
(\al+G)\cap R=(\al+G)\cap R^\prime=(\al+G)\cap
R(k,0)=\al+\BBbb{Z}k(\al)a
\end{equation} for $\al\in\Pi$.
To complete the proof, it suffices to show that
\begin{equation}\nonumber
\forall \al\in\Pi,\,R^{\prime\prime}\cap(2\al +G)\subset
2\al
+\BBbb{Z}k(\al)a\,.
\end{equation}
Let $\beta\in R^{\prime\prime}\cap(2\al +G)$.
Let $x:=\beta-2\al$.
Then
$\sigma_\al\sigma_\beta(\al)=3\al-\beta=\al-x$. By
(\ref{eqn:alGcapR}),
we have $x\in \BBbb{Z}k(\al)a$, as desired. \hfill $\BBox$

\section{Elliptic Lie algebras and superalgebras}

\subsection{Definition with generators and relations}

Let ${\cal D}$ be a QEBS
with $l\geq 2$. For $\al\in\Pi$, let
\begin{equation}
c(\al ):=\left\{\begin{array}{ll}
2 & \mbox{if $g(\al)=\BBbb{Z}$ or $2\BBbb{Z}+1$,} \\
1 & \mbox{otherwise,} \\
\end{array}\right.
\nonumber
\end{equation} and let
\begin{equation}
\al^*:=c(\al)\al+k(\al)a\in{\cal E}.
\nonumber
\end{equation}
For $\al\in\Rkg$, let
\begin{equation}
p(\al ):=\left\{\begin{array}{ll}
1 & \mbox{if $2\al\in\Rkg$,} \\
0 & \mbox{otherwise.} \\
\end{array}\right.
\nonumber
\end{equation}
Let
\begin{equation}
{\cal
A}:=\{(\al,\beta,y)\in\Pi\times\Pi\times\BBbb{N}\,|\,\al\ne\beta,\,
J(\al,\beta^\vee)=-1,\,k(\al)y=k(\beta)\}.
\nonumber
\end{equation}\par
For a subset $T$ of the ${\cal E}$, let $-T:=\{-t\in{\cal
E}|t\in T\}$.
Let $\Pi^*:=\{\al^*|\al\in\Pi\}$ and
${\cal B}_+:=\Pi\cup\Pi^*$. Let
\begin{equation}
{\cal B}:=
{\cal B}_+\cup -{\cal B}_+.
\nonumber
\end{equation}
Let $({\cal B}\times{\cal B})^\prime:=
\{(\mu,\nu)\in {\cal B}\times{\cal
B}|\mu\ne\nu,\,\mu+\nu\ne 0\}$.
For $(\mu,\nu)\in({\cal B}\times{\cal B})^\prime$, let
\begin{equation}
x_{\mu,\nu}:=
\left\{\begin{array}{ll}
1-J(\mu^\vee,\nu)  &  \mbox{if
$J(\mu^\vee,\nu)<0$,} \\ 0  &  \mbox{if
$J(\mu^\vee,\nu)\geq 0$.}
\end{array}\right.
\nonumber
\end{equation}

Keep the notation as above. Let ${\frak g}_{{\cal
D}}={\frak g}_\EPiakg$ be
the Lie superalgebra defined by
generators
\begin{equation}\label{eqn:generators}
h_\sigma\,(\sigma\in{\cal E}), E_\mu \,(\mu\in{\cal B})
\end{equation} with parities
\begin{equation}\label{eqn:generatorsparities}
p(h_\sigma)=0,\,p(E_\mu)=p(\mu)
\end{equation} and the following defining relations.
\newline\par
({\rm SR}1) \quad
$xh_\sigma+yh_\tau=h_{x\sigma+y\tau}$ \quad
if $x,\,y\in\BBbb{C}$ and $\sigma,\,\tau\in{\cal E}$, \par
({\rm SR}2) \quad
     $[h_\sigma,h_\tau]=0$ \quad if $\sigma,\,\tau\in{\cal
E}$, \par
({\rm SR}3) \quad
$[h_\sigma,E_\mu]=J(\sigma,\mu)E_\mu$ \quad
if $\sigma\in{\cal E}$ and $\mu\in{\cal B}$, \par
({\rm SR}4) \quad
$[E_\mu,E_{-\mu}]=h_{\mu^\vee}$ \quad
if $\mu\in{\cal B}_+$, \par
({\rm SR}5) \quad
$({\rm ad}E_\mu)^{x_{\mu,\nu}}E_\nu=0$ \quad
if $(\mu,\nu)\in( {\cal B}\times{\cal B})^\prime$, \par
({\rm SR}6) \quad
$c(\al)({\rm ad}E_{\al^*})^yE_\beta=({\rm
ad}E_\al)^{c(\al)y}E_{\beta^*}$
\quad
if $(\al,\beta,y)\in{\cal A}$, \par
({\rm SR}7) \quad
$(-1)^{c(\al)+1}c(\al)({\rm ad}E_{-\al^*})^yE_{-\beta}=
({\rm ad}E_{-\al})^{c(\al)y}E_{-\beta^*}$
\quad
if $(\al,\beta,y)\in{\cal A}$, \par
({\rm SR}8) \quad
$({\rm ad}E_\al)^i({\rm ad}E_{\al^*})^{y-i}E_\beta=0$
\quad
if $(\al,\beta,y)\in{\cal A}$ and $1\leq i\leq y-1$, \par
({\rm SR}9) \quad
$({\rm ad}E_{-\al})^i({\rm
ad}E_{-\al^*})^{y-i}E_{-\beta}=0$ \quad
if $(\al,\beta,y)\in{\cal A}$ and $1\leq i\leq y-1$.
\newline\par
We call the ${\frak g}_{{\cal D}}$ the {\it elliptic Lie
(super)algebra}. In Introduction, the ${\frak g}_{{\cal
D}}$
is also denoted as ${\frak g}_R$, where $R=\Rkg$.
For $\mu\in{\cal E}$, let
${\frak g}_{{\cal D},\mu}:=\{X\in{\frak g}_{{\cal D}}|
\forall h_\sigma,\, [h_\sigma,X]=J(\sigma,\mu)X\}$.
Define
the sub Lie superalgebra ${\frak h}_{{\cal D}}$ of
${\frak g}_{{\cal D}}$ by ${\frak h}_{{\cal D}}:=
\{h_\sigma\in{\frak g}_{{\cal D},0}|\sigma\in{\cal E}\}$.
\begin{lemma}\label{lemma:Cartannondeg}
Keep the notation as above.
Then $h_\sigma\ne 0$ for $\sigma\in{\cal
E}\setminus\{0\}$.
In particular, $\dim{\frak h}_{{\cal D}}=l+4$.
\end{lemma}
\par Proof of the lemma shall be given in Subsec.~4.2.

\subsection{Main theorem}

We see that there exists
a unique $\delta\in\BBbb{Z}_+\Pi$
such that $\BBbb{Z}\delta
=\{\lambda\in\BBbb{Z}\Pi|J(\lambda,\lambda)=0\}$.
By (\ref{eqn:lambdazeroalzero}),
we have
$J(\Lambda_\delta,\delta)=1$. \par
Let $\BBbb{Z}^{2,\prime}:=\BBbb{Z}^2\setminus\{(0,0)\}$.

\begin{theorem}\label{theorem:main}  Let ${\cal
D}=\EPiakg$ be a QEBS with
$l\geq 2$.
Then we have
\begin{equation}\nonumber
{\frak g}_{{\cal D}}={\frak h}_{{\cal D}}\bigoplus
\Bigl(\bigoplus_{\nu\in R(k,g)}{\frak g}_{{\cal
D},\nu}\Bigr)\bigoplus
\Bigl(\bigoplus_{(m,n)\in \BBbb{Z}^{2,\prime}}{\frak
g}_{{\cal
D},m\delta+na}\Bigr).
\end{equation} Moreover ${\frak h}_{{\cal D}}={\frak
g}_{{\cal D},0}$, and
$\dim {\frak g}_{{\cal D},\al}=1$ for all $\al\in R(k,g)$.
\end{theorem}
\par Proof of the theorem shall be given in Subsec.~3.2.
\newline\par From now until the end of Subsec.~2.2, we
suppose that we have
proved
Theorem~\ref{lemma:Cartannondeg}. Let $\nu\in R(k,g)$.
Let $E^\prime_{\pm\nu}\in{\frak g}_{{\cal D},\pm\nu}$
be such that
$[E^\prime_\nu,E^\prime_{-\nu}]=h_{\nu^\vee}$.
By Theorem~\ref{theorem:main}, $E^\prime_{\pm\nu}$ are
locally nilpotent
since
$|J(\sigma+r\nu,\sigma+r\nu))|\rightarrow +\infty$ as
$r\rightarrow +\infty$, where $\sigma\in{\cal E}$.
Hence we can define $n_\nu =n_{E^\prime_\nu}
\in{\rm Aut}({\frak g}_{{\cal D}})$
by
\begin{equation}\label{eqn:nmu}
n_\nu
=\left\{\begin{array}{ll}
\exp{\rm ad}E^\prime_\nu\exp{\rm
ad}(-E^\prime_{-\nu})\exp{\rm
ad}E^\prime_\nu
& \mbox{if $p(\nu)=0$,} \\
\exp ({\frac 1 4}{\rm ad}[E^\prime_\nu,E^\prime_\nu])\exp
({\frac 1 4}{\rm
ad}[E^\prime_{-\nu},E^\prime_{-\nu}])
\exp ({\frac 1 4}{\rm ad}[E^\prime_\nu,E^\prime_\nu])
& \mbox{if $p(\nu)=1$.}
\end{array}\right.
\end{equation}

\begin{theorem}\label{theorem:submain} Let ${\cal D}$ be a
PEBS with $l\geq
2$. Then ${\cal D}$ is an EBS if and only if it is a QEBS.
\end{theorem}
      {\it Proof.} The `only-if'-part follows from
Lemma~\ref{lemma:kg}.
Here we prove the `if'-part. Recall $n_\nu$
($\nu\in\Rkg$) from (\ref{eqn:nmu}).
Notice that $n_\nu({\frak g}_{{\cal D},\lambda})
={\frak g}_{{\cal D},s_\nu(\lambda)}$. Then, by
Theorem~\ref{theorem:main},
we see that
${\cal D}$ is an EBS.
\hfill $\BBox$

\section{Proof of Theorem~\ref{theorem:main}}
In this section, we suppose that we have proved
Lemma~\ref{lemma:Cartannondeg}.

\subsection{Rank one and two subsystems}
Let $S$ be a finite subset of $\Rkg$.
Assume that the elements of the $S$ are linearly
independent and that the square matrix
$A_S:=(J(\al^\vee,\beta))_{\al,\beta\in S}$
is an affine type generalized Cartan matrix in the sense
of \cite[\S
4.8]{kac1}.
Then we call the $S$ the {\it affine type subset} of
$\Rkg$.
Let $EW_S$ be the subgroup of ${\rm GL}({\cal E})$
generated by $s_\mu$ ($\mu\in S$).
Let $S^{{\rm odd}}:=\{\al\in S|p(\al)=1\}$
and let
\begin{equation}
      \Rkg^S:=\bigcup_{w\in EW_S}w\Bigl(\bigcup_{\al\in
S\setminus S^{{\rm
odd}}}\{\al\}\bigcup\bigcup_{\al\in S^{{\rm
odd}}}\{\al,2\al\}\Bigr).
\nonumber
\end{equation} Then we see that $\Rkg^S$ is the affine
type
(real) root system with the base $S$; $\Rkg^S$ is reduced
if and only if $S^{{\rm odd}}=\emptyset$.
For the pair $(A_S,S^{{\rm odd}})$ of the above $A_S$ and
$S^{{\rm odd}}$,
we define the Dynkin diagram
$\Gamma(A_S,S^{{\rm odd}})$
in the same manner as in \cite{kac2}. If
$\Gamma(A_S,S^{{\rm odd}})$ is
called ${\bf X}$ in the tables \cite[Tables 1-4]{kac2},
we say that the name of $S$ is ${\bf X}$.
The following two lemmas follow from ({\rm KG}1-3),
Lemma~\ref{lemma:Rkgs}
and the well-known fact \cite[Appendixes 1-2]{MacDonald}.

\begin{lemma}\label{lemma:rankone} Let ${\cal D}$ be a
QEBS
with $l\geq 2$. Let $\al\in\Pi$.
Then $\{\al,-\al^*\}$ is an affine type subset of $\Rkg$,
and we have
$\Rkg_{\{\al\}}=\Rkg^{\{\al,-\al^*\}}$.
Moreover, letting ${\bf X}$ be the name of
the $\{\al,-\al^*\}$, we have one of the following cases.
\par
(i) $g(\al)=\emptyset$, ${\bf X}={\bf A}_1^{(1)}$. \par
(ii) $g(\al)=2\BBbb{Z}+1$, ${\bf X}={\bf A}_2^{(2)}$. \par
(iii) $g(\al)=\BBbb{Z}$, ${\bf X}={\bf B}^{(1)}(0,1)$.
\par
(iv) $g(\al)=2\BBbb{Z}$, ${\bf X}={\bf C}^{(2)}(2)$. \par
(v) $g(\al)=4\BBbb{Z}+2$, ${\bf X}={\bf A}^{(4)}(0,2)$,
$p(\al)=0$. \par
(vi) $g(\al)=4\BBbb{Z}$, ${\bf X}={\bf A}^{(4)}(0,2)$,
$p(\al)=1$.
\end{lemma}

\begin{lemma}\label{lemma:ranktwo} Let ${\cal D}$ be a
QEBS with $l\geq 2$. Let $\al$, $\beta\in\Pi$
be such that $J(\al,\beta^\vee)=-1$. Then
$g(\beta)=\emptyset$
and there exists a unique $\gamma\in\Rkg_{\{\al,\beta\}}
\cap(-a+\BBbb{Z}_-\Pi)$ such that
$\{\al,\beta,\gamma\}$ is an affine type subset of the
$\Rkg_{\{\al,\beta\}}$ Moreover we have
$\Rkg^{\{\al,\beta,\gamma\}}
=\Rkg_{\{\al,\beta\}}$. Furthermore, letting
${\bf Y}$ be the name of the $\{\al,\beta,\gamma\}$
and letting $Jkg:=\{J(\al^\vee,\beta),{\frac
{k(\beta)}{k(\al)}},
g(\al)\}$,
we have one of the following cases.
\par
(i) $Jkg=\{-1,1,\emptyset\}$,
$\gamma=s_\al(-\beta^*)$,
${\bf Y}={\bf A}_2^{(1)}$. \par
(ii) $Jkg=\{-2,1,\emptyset\}$,
$\gamma=s_\al(-\beta^*)$,
${\bf Y}={\bf C}_2^{(1)}$. \par
(iii) $Jkg=\{-3,1,\emptyset\}$,
$\gamma=s_\beta s_\al(-\beta^*)$,
${\bf Y}={\bf G}_2^{(1)}$. \par
(iv) $Jkg=\{-2,2,\emptyset\}$,
$\gamma=s_\beta(-\al^*)$,
${\bf Y}={\bf D}_3^{(2)}$. \par
(v) $Jkg=\{-3,3,\emptyset\}$,
$\gamma=s_\al s_\beta(-\al^*)$,
${\bf Y}={\bf D}_4^{(3)}$. \par
(vi) $Jkg=\{-2,1,2\BBbb{Z}+1\}$,
$\gamma=s_\beta(-\al^*)$,
${\bf Y}={\bf A}_4^{(2)}$. \par
(vii) $Jkg=\{-2,1,\BBbb{Z}\}$,
$\gamma=s_\beta(-\al^*)$,
${\bf Y}={\bf B}^{(1)}(0,2)$. \par
(viii) $Jkg=\{-2,1,2\BBbb{Z}\}$,
$\gamma=s_\al(-\beta^*)$,
${\bf Y}={\bf A}^{(2)}(0,3)$. \par
(ix) $Jkg=\{-2,2,2\BBbb{Z}\}$,
$\gamma=s_\beta(-\al^*)$,
${\bf Y}={\bf C}^{(2)}(3)$. \par
(x) $Jkg=\{-2,2,4\BBbb{Z}+2\}$,
$\gamma=s_\beta(-\al^*)$,
${\bf Y}={\bf A}^{(4)}(0,4)$. \par
(xi) $Jkg=\{-2,2,4\BBbb{Z}\}$,
$\gamma=s_\beta(-\al^*)$,
${\bf Y}={\bf A}^{(4)}(0,4)$.
\end{lemma}

\subsection{Embedded rank two affine (super)algebras}
Let ${\cal D}$ be a QEBS with $l\geq 2$.
Let $\mu$, $\nu\in{\cal B}_+$ be such that
$J(\mu^\vee,\nu)=-1$.
Then we have
\begin{equation}\label{eqn:valueofn}
\left\{
\begin{array}{l}
n_{xE_\mu}(E_{\pm \nu})=\mp x^{\pm 1}[E_{\pm
\mu},E_{\pm \nu}],\\
n_{xE_\nu}(E_{\pm \mu})=(\pm x)^{\mp J(\nu^\vee,\mu)}
{\frac {1} {(-J(\nu^\vee,\mu))!}}
({\rm ad}E_{\pm \nu})^{-J(\nu^\vee,\mu)}E_{\pm \mu}
\end{array}\right.
\end{equation} for $x\in\BBbb{C}\setminus\{0\}$ (see
(\ref{eqn:nmu})). \par
    From now until the end of the paper, we assume that
for every $\mu\in{\cal B}$,
$n_\mu$ always denotes $n_{E_\mu}\in{\rm Aut}({\frak
g}_{{\cal D}})$
and we let $X\sim Y$ mean that there exists a
$z\in\BBbb{C}\setminus\{0\}$
such that
$X=zY$.
\begin{lemma}\label{lemma:embedinglemma}
Keep the notation as in Lemma~\ref{lemma:ranktwo}.
Let $T:=\{\al,\beta,\gamma\}$.
Assume $\gamma=s_{\nu_1}\cdots s_{\nu_r}(-\mu^*)$ to be
the same
expression as
in Lemma~\ref{lemma:ranktwo} (i)-(xi). Set
\begin{equation}\label{eqn:twoserre}
E_{\pm\gamma}:=(\pm 1)^{p(\mu^*)}
n_{\nu_1}\cdots n_{\nu_r}(E_{\mp\mu^*}).
\end{equation}
Then we have
\begin{equation}\label{eqn:ranktwoserre}
\left\{\begin{array}{ll}
({\rm ad}E_{\pm\mu})^{1-J(\mu^\vee,\nu)}E_{\pm\nu}=0
& \mbox{if $\mu$, $\nu\in T$ with $\mu\ne\nu$,} \\
\mbox{$[E_\mu,E_{-\nu}]=\delta_{\mu\nu}h_{\mu^\vee}$}
&
      \mbox{if $\mu$, $\nu\in T$.}
\end{array}\right.
\end{equation}
\end{lemma}\par
{\it Proof.} If the name of $T$ is neither ${\bf
A}_4^{(2)}$ nor
${\bf B}^{(1)}(0,2)$, the equalities
(\ref{eqn:ranktwoserre})
is proved in a similar way to \cite[\S 2.3]{yamane2}.
We assume that the name of $T$ is ${\bf A}_4^{(2)}$ or
${\bf
B}(0,2)^{(1)}$.
It is clear that
$[E_\gamma,E_{-\gamma}]=h_{s_\beta(-(\al^*)^\vee)}=h_{\gamma^\vee}$.
We have
\begin{equation}\nonumber
E_{\gamma} = n_\beta(E_{-\al^*})
      \sim [E_{-\beta},[E_{-\beta},E_{-\al^*}]]
      \sim
[E_{-\beta},[E_{-\al},[E_{-\al},E_{-\beta^*}]]].
\end{equation}
It is clear that $[E_\gamma,E_{-\beta}]=0$. We have
\begin{eqnarray*}
[E_\gamma,E_{-\al}] &\sim & [[E_{-\beta},E_{-\al}],
[E_{-\al},[E_{-\al},E_{-\beta^*}]]] \\
&\sim & [[E_{-\beta},E_{-\al}],[E_{-\beta},E_{-\al^*}]] \\
&\sim & n_\beta([E_{-\al},[E_{-\beta},E_{-\al^*}]]) \\
&\sim &
n_\beta([[E_{-\al},[E_{-\al},[E_{-\al},E_{-\beta^*}]]])=
0.
\end{eqnarray*}
Similarly we have $[E_{-\gamma},E_\beta]=0$ and
$[E_{-\gamma},E_\al]=0$.
Since the $E_\lambda$ ($\lambda\in T\cup -T$)
are locally nilpotent,
we have all the equalities in (\ref{eqn:ranktwoserre}).
\hfill $\BBox$ \newline\par

For a subset $S$ of $\Pi$, let ${\frak g}_{{\cal D}}^S$
be the sub Lie superalgebra of ${\frak g}_{{\cal D}}$
generated by $h_\sigma$ ($\sigma\in{\cal E}$)
and $E_\mu$, $E_{-\mu}$, $E_{\mu^*}$, $E_{-\mu^*}$
($\mu\in S$). For $\nu\in{\cal E}$, let ${\frak g}_{{\cal
D},\nu}^S
:={\frak g}_{{\cal D}}^S\cap{\frak g}_{{\cal D},\nu}$.
Let $m_S:=\min\{k(\al)|\al\in S\}$.

\begin{lemma}\label{lemma:sroots}
Let ${\cal D}$ be a QEBS with $l\geq 2$.
Let $S$ be a subset of $\Pi$ such that $|S|=1$ or $2$.
Then we have
\begin{equation}\label{eqn:rootdecompoitionofs}
{\frak g}_{{\cal D}}^S={\frak h}_{{\cal D}}\bigoplus
\Bigl(\bigoplus_{\lambda\in R(k,g)_S}{\frak g}_{{\cal
D},\lambda}^S\Bigr)\bigoplus
\Bigl(\bigoplus_{n\in \BBbb{Z}
\setminus\{0\}}{\frak g}_{{\cal D},nm_Sa}^S\Bigr).
\end{equation}
Moreover we have $\dim{\frak g}_{{\cal D},\lambda}^S=1$
for $\lambda\in R(k,g)_S$. Furthermore
\begin{equation}\label{eqn:gdnmsaalgdnmsabeta}
{\frak g}_{{\cal D},nm_Sa}^S={\frak g}_{{\cal
D},nm_Sa}^{\{\al\}}
+{\frak g}_{{\cal D},nm_Sa}^{\{\beta\}}
\end{equation} if $S=\{\al,\beta\}$.
\end{lemma}\par
{\it Proof.}
We first assume $|S|=1$ and $S=\{\al\}$. Recall the affine
subset $U:=\{\al,-\al^*\}$ from Lemma~\ref{lemma:rankone}.
It follows from Lemma~\ref{lemma:Cartannondeg}
that $E_{\pm\mu}\ne 0$ for $\mu\in U$, since
$[E_\mu,E_{-\mu}]
=h_{\mu^\vee}\ne 0$.
Since
$E_{\pm\mu}$ ($\mu\in U$) satisfy the Serre relations
(see ({\rm SR}5)), the lemma
follows from the well-known argument in the proof of
\cite[Corollary~5.12]{kac1}
(see also \cite[Proposition~1.6]{kac2}). \par
Assume $|S|=2$, and $S=\{\al,\beta\}$.
If $J(\beta^\vee,\al)=0$, the lemma follows from the same
argument as
above.
Assume $J(\beta^\vee,\al)=-1$.
Recall $\gamma$ and $E_{\pm\gamma}$ from
Lemma~\ref{lemma:embedinglemma}; especially recall
(\ref{eqn:twoserre}). Let $T:=\{\al,\beta,\gamma\}$.
Notice that
$E_{\pm\gamma}\in{\frak g}_{{\cal D}}^S$.
Let ${\frak g}_{{\cal D}}^{(T)}$ be the sub Lie
superalgebra
of ${\frak g}_{{\cal D}}^S$ generated by ${\frak h}_{{\cal
D}}$
and $E_{\pm\omega}$ ($\omega\in T$). Then
$E_{\pm\mu^*}=n_{\nu_r}^{-1}\cdots
n_{\nu_1}^{-1}E_{\mp\gamma}\in{\frak g}_{{\cal
D}}^{(T)}$.
Let $\rho\in\{\al,\beta\}$ be such that
$\rho\ne\mu$.
By ({\rm SR}6-7)
and by (\ref{eqn:valueofn}),
we have $n_\mu E_{\pm\rho^*}=n_{\mu^*} E_{\pm\rho}$.
Hence $E_{\pm\rho^*}=n_\mu^{-1}n_{\mu^*}
E_{\pm\rho}\in{\frak g}_{{\cal
D}}^{(T)}$.
Hence ${\frak g}_{{\cal D}}^{(T)}={\frak g}_{{\cal D}}^S$.
By Lemma~\ref{lemma:embedinglemma}, $E_{\pm\omega}$
($\omega\in T$)
satisfy the Serre relations (\ref{eqn:ranktwoserre}).
Using the well-known argument in the proof of
\cite[Corollary~5.12]{kac1}
again,
we have (\ref{eqn:rootdecompoitionofs}) for the $S$.
The equality (\ref{eqn:gdnmsaalgdnmsabeta})
follows
from the fact that ${\frak g}_{{\cal D}}^S$ is
generated by ${\frak g}_{{\cal D},\lambda}^{\{\al\}}$
($\lambda\in
R(k,g)_{\{\al\}}$)
and ${\frak g}_{{\cal D},\nu}^{\{\beta\}}$ ($\nu\in
R(k,g)_{\{\beta\}}$)
(See also Lemma~\ref{lemma:Rkgs}).
This completes the proof.
\hfill $\BBox$ \newline\par
Keep the notation as above. Let $S$ be a subset of $\Pi$.
Let $ED_S:=\BBbb{Z}S+\BBbb{Z}m_Sa$.
Define the subsets $ED_{S,+}$ and $ED_{S,-}$ of the $ED_S$
by $ED_{S,\pm}:=(\BBbb{Z}_\pm
S+\BBbb{Z}m_Sa)\setminus\BBbb{Z}m_Sa$.
Define the sub Lie superalgebras
${\frak n}_{{\cal D}}^{S,+}$, ${\frak n}_{{\cal D}}^{S,-}$
${\frak l}_{{\cal D}}^{S,+}$ and ${\frak l}_{{\cal
D}}^{S,-}$ of the
${\frak g}_{{\cal D}}^S$
by
\begin{equation}\nonumber
{\frak n}_{{\cal D}}^{S,\pm}:=\bigoplus_{\lambda\in
ED_{S,\pm}}
{\frak g}_{{\cal D},\lambda}^S
\end{equation} and \begin{equation}\nonumber
{\frak l}_{{\cal D}}^{S,\pm}:=\bigoplus_{n\in
\BBbb{Z}_\pm\setminus\{0\}}
{\frak g}_{{\cal D},nm_sa}^S.
\end{equation}
\begin{lemma}\label{lemma:striangular}
Let $S$ be a subset of $\Pi$. Then the following hold.
\par
(1) ${\frak g}_{{\cal D}}^S={\frak n}_{{\cal D}}^{S,+}
\oplus{\frak l}_{{\cal D}}^{S,+}\oplus{\frak h}_{{\cal D}}
\oplus{\frak l}_{{\cal D}}^{S,-}\oplus{\frak n}_{{\cal
D}}^{S,-}$.\par
(2) The ${\frak n}_{{\cal D}}^{S,+}$
is generated by
      ${\frak n}_{{\cal D}}^{\{\al\},+}$
with $\al\in S$. The ${\frak n}_{{\cal D}}^{S,-}$
is generated by
      ${\frak n}_{{\cal D}}^{\{\al\},-}$
with $\al\in S$.
\par
(3) ${\frak l}_{{\cal D}}^{S,+}
=\sum_{\al\in S}{\frak l}_{{\cal D}}^{\{\al\},+}$,
${\frak l}_{{\cal D}}^{S,-}
=\sum_{\al\in S}{\frak l}_{{\cal D}}^{\{\al\},-}$.
\par
(4) $\dim{\frak g}_{{\cal D},\lambda}^S=1$ for $\lambda\in
R(k,g)_{\{\al\}}$
if $\al\in\Pi$.
\end{lemma} \par {\it Proof.}
If $|S|=1$ or $2$, the lemma follows from
Lemma~\ref{lemma:sroots}.
Assume $|S|\geq 3$.
Notice that the ${\frak g}_{{\cal D}}^S$ is generated
by ${\frak h}_{{\cal D}}$ and ${\frak n}_{{\cal
D}}^{\{\gamma\},+}$,
${\frak n}_{{\cal D}}^{\{\gamma\},-}$ with $\gamma\in S$.
Assume $\al$, $\beta\in S$ with $\al\ne\beta$.
Notice that
\begin{equation}\label{eqn:zalzbetaza}
(\BBbb{Z}_\pm\al+\BBbb{Z}_\mp\beta+\BBbb{Z}a)\cap
R(k,g)_S=\emptyset.
\end{equation}
    From Lemma~\ref{lemma:sroots} and
(\ref{eqn:zalzbetaza})
when
$J(\al^\vee,\beta)\ne 0$,
or from the defining relations ({\rm SR}1-9) when
$J(\al^\vee,\beta)=0$,
it follows that
$[{\frak n}_{{\cal D}}^{\{\al\},+},{\frak n}_{{\cal
D}}^{\{\beta\},-}]=\{0\}$ and
$[{\frak n}_{{\cal D}}^{\{\beta\},\pm},
{\frak l}_{{\cal D}}^{\{\al\},+}+{\frak l}_{{\cal
D}}^{\{\al\},-}]
\subset{\frak n}_{{\cal D}}^{\{\beta\},\pm}$.
Then the lemma follows from this fact and
Lemmas~\ref{lemma:Cartannondeg} and \ref{lemma:sroots}.
\hfill
$\BBox$
\newline\par
{\it Proof of Theorem~\ref{theorem:main}.} The theorem
follows from
Lemmas~\ref{lemma:zsetimes},
\ref{lemma:Rkgs} and \ref{lemma:striangular}
and from the existence of the $n_\gamma$'s
with $\gamma\in\Pi$. \hfill $\BBox$

\section{Proof of Lemma~\ref{lemma:Cartannondeg}}

\subsection{Contragredient Lie superalgebra}

Here we first recall the definition of the
contragredient Lie superalgebras \cite{kac2}. Let ${\bar
I}$
be a finite set.
Let ${\bar I}^{{\rm odd}}$ be a subset of ${\bar I}$.
Define a map ${\bar p}:{\bar I}\rightarrow\{0,1\}$ by
${\bar p}(i)=1$ ($i\in{\bar I}^{{\rm odd}}$) and
${\bar p}(j)=0$ ($j\in{\bar I}\setminus{\bar I}^{{\rm
odd}}$).
Let ${\bar A}:=({\bar a}_{ij})_{i,j\in{\bar I}}$
be an ${\bar I}\times{\bar I}$ matrix.
Let ${\bar {\frak H}}$ be the $2|{\bar I}|$-dimensional
$\BBbb{C}$-vector space.
Let ${\bar {\frak H}}^*$ be the dual space of ${\bar
{\frak H}}$.
Let $\{{\bar \al}_i,\,{\bar \gamma}_i\,(i\in{\bar I})\}$
be a basis of ${\bar {\frak H}}$.
Let $\{{\bar h}_i,\,{\bar t}_i\,(i\in{\bar I})\}$
be a basis of ${\bar {\frak H}}^*$.
We assume that ${\bar \al}_i({\bar h}_j)={\bar a}_{ji}$,
${\bar \al}_i({\bar t}_j)=
{\bar \gamma}_i({\bar h}_j)=\delta_{ji}$ and
${\bar \gamma}_i({\bar t}_j)=0$.
Let  ${\bar \Pi}:=\{{\bar \al}_i\,(i\in{\bar I})\}$.
Let ${\bar \Pi}^\vee:=\{{\bar h}_i\,(i\in{\bar I})\}$.
For the datum
${\bar {\frak D}}:=({\bar A},{\bar I},{\bar I}^{{\rm
odd}},{\bar
\Pi},{\bar
\Pi}^\vee)$, we define a Lie superalgebra
${\bar {\frak G}}^\prime_{{\bar {\frak D}}}:={\bar {\frak
G}}^\prime({\bar A},{\bar
I}^{{\rm
odd}})$ by generators
\begin{equation}
{\bar h}^\prime_i,\,{\bar t}^\prime_i,\
{\bar E}^\prime_i,\,{\bar F}^\prime_i\,(i\in{\bar I})
\nonumber
\end{equation} with parities
\begin{equation}
p({\bar h}^\prime_i)=p({\bar t}^\prime_i)=0,\,\,p({\bar
E}^\prime_i)=p({\bar
F}^\prime_i)={\bar p}(i)
\nonumber
\end{equation} and defining relations
\begin{equation}
\left\{\begin{array}{l}
\mbox{$[{\bar h}^\prime_i,{\bar h}^\prime_j]=
[{\bar h}^\prime_i,{\bar t}^\prime_j]=
[{\bar t}^\prime_i,{\bar t}^\prime_j]=0$,} \\
\mbox{$[{\bar h}^\prime_i,{\bar E}^\prime_j]={\bar
a}_{ij}{\bar
E}^\prime_j$,\,$[{\bar
t}^\prime_i,{\bar E}^\prime_j]=\delta_{ij}{\bar
E}^\prime_j$ ,\,$[{\bar
h}^\prime_i,{\bar F}^\prime_j]=-{\bar a}_{ij}{\bar
F}^\prime_j$ ,\,$[{\bar
t}^\prime_i,{\bar F}^\prime_j]=-\delta_{ij}{\bar
F}^\prime_j$,}  \\
\mbox{$[{\bar E}^\prime_i,{\bar
F}^\prime_j]=\delta_{ij}{\bar
h}^\prime_i$\,.}
\end{array}\right.
\nonumber
\end{equation}
Let ${\bar {\frak H}}^\prime$ be the sub Lie superalgebra
of
${\bar {\frak G}}^\prime_{{\bar {\frak D}}}$
generated by ${\bar h}^\prime_i$, ${\bar t}^\prime_i$.
Let ${\frak r}$ be the ideal of ${\bar {\frak
G}}^\prime_{{\bar {\frak D}}}$ which
is
maximal
among the ones ${\frak r}^\prime$ such that ${\frak
r}^\prime\cap{\bar
{\frak H}}^\prime=\{0\}$. We denote by ${\bar {{\frak
G}}}_{{\bar {\frak D}}}={\bar
{{\frak G}}}({\bar A},{\bar I}^{{\rm
odd}})$ the quotient Lie superalgebra
${\bar {\frak G}}^\prime_{{\bar {\frak D}}}/{\frak r}$.
In this paper,
we call the ${\bar {{\frak G}}}_{{\bar {\frak D}}}$
the {\it contragredient
Lie superalgebra}. Let ${\bar \pi}:{\bar {\frak
G}}^\prime_{{\bar {\frak D}}}\rightarrow
{\bar {{\frak G}}}_{{\bar {\frak D}}}$ be a natural
projective
map.
Notice that $\dim{\bar \pi}({\bar {\frak
H}}^\prime)=2|{\bar I}|$.
By abuse of notation, we shall also denote
${\bar \pi}({\bar {\frak H}}^\prime)$,
${\bar \pi}({\bar h}^\prime_i)$ and ${\bar \pi}({\bar
t}^\prime_i)$
by
${\bar {\frak H}}$,
${\bar h}_i$ and ${\bar t}_i$,
respectively.
We shall also denote
${\bar \pi}({\bar E}^\prime_i)$ and ${\bar \pi}({\bar
F}^\prime_i)$
by
${\bar E}_i$ and ${\bar F}_i$,
respectively.
\par
Keep the notation as above.  Let ${\bar I}^{{\rm
pos}}:=\{i\in{\bar
I}|{\bar
a}_{ii}\ne 0\}$. Let ${\bar I}^{{\rm null}}:={\bar
I}\setminus{\bar
I}^{{\rm pos}}$.
Define the square matrix
${\bar A}^{{\rm pos}}$ by ${\bar A}^{{\rm pos}}:=({\bar
a}_{ij})_{i,j\in{\bar
I}^{{\rm pos}}}$. We say that ${\bar {\frak D}}$ is a {\it
handy
datum} if
the following hold. \par
({\rm HD}1) If $i\in{\bar I}^{{\rm pos}}$, then ${\bar
a}_{ii}=2$.\par
({\rm HD}2) If $i$, $j\in{\bar I}^{{\rm pos}}$ with $i\ne
j$
and $0\leq |{\bar a}_{ij}|\leq |{\bar a}_{ij}|$,
then $({\bar a}_{ij},{\bar a}_{ji})$ is
$(0,0)$,
$(-1,-1)$, $(-1,-2)$ or $(-1,-3)$. \par
({\rm HD}3) If $i$, $j\in{\bar I}^{{\rm null}}$ with $i\ne
j$,
then $({\bar a}_{ij},{\bar a}_{ji})=(0,0)$ or $(2,2)$.
\par
({\rm HD}4) If $i\in{\bar I}^{{\rm pos}}$ and $j\in{\bar
I}^{{\rm null}}$,
then $({\bar a}_{ij},{\bar a}_{ji})=(0,0)$, $(-1,-1)$ or
$(-1,-2)$.  \par
({\rm HD}5)
If $i\in{\bar I}^{{\rm pos}}$ and $j\in{\bar I}^{{\rm
null}}$,
then  $({\bar a}_{ij},{\bar a}_{ji})=(-1,-1)$
if and only if there exists an $r\in{\bar I}^{{\rm
null}}\setminus\{j\}$
such that ${\bar a}_{ir}\ne 0$ and ${\bar a}_{jr}\ne 0$.
\par
({\rm HD}6) ${\bar I}^{{\rm null}}\subset{\bar I}^{{\rm
odd}}$. \par
({\rm HD}7) If $i\in{\bar I}^{{\rm pos}}\cap{\bar I}^{{\rm
odd}}$ and $j\in{\bar I}^{{\rm null}}$, then
$({\bar a}_{ij},{\bar a}_{ji})=(0,0)$.
\par ({\rm HD}8) If $i\in{\bar I}^{{\rm pos}}\cap{\bar
I}^{{\rm
odd}}$, $j\in{\bar I}^{{\rm pos}}\setminus\{i\}$ and
${\bar a}_{ij}\ne 0$,
then $j\notin{\bar I}^{{\rm odd}}$ and
$({\bar a}_{ij},{\bar a}_{ji})=(-2,-1)$.
\par ({\rm HD}9) If $i\in{\bar I}^{{\rm null}}$, then
there exists a unique
$j\in{\bar I}^{{\rm null}}$ such that $i\ne j$ and ${\bar
a}_{ij}\ne 0$.
\par ({\rm HD}10) There exist ${\bar
\epsilon}_i\in\BBbb{C}\setminus\{0\}$
($i\in{\bar I}$) such that
${{\bar D}}^{-1}{{\bar A}}$ is a symmetric matrix,
where ${{\bar D}}$ is the diagonal matrix
$(\delta_{ij}{\bar \epsilon}_i)$.
\newline\par
Assume  ${\bar {\frak D}}$ to be a handy datum. Then
there exists a nondegenerate symmetric bilinear form
${\bar J}:{\bar {\frak H}}^*\times{\bar {\frak
H}}^*\rightarrow\BBbb{C}$
such that ${\bar J}({\bar \al}_i,{\bar \al}_j)={\bar
\epsilon}_i^{\,-1}{\bar
a}_{ij}$,
${\bar J}({\bar \gamma}_i,{\bar \al}_j)={\bar
\epsilon}_i^{\,-1}\delta_{ij}$,
${\bar J}({\bar \gamma}_i,{\bar \gamma}_j)=0$
For $\sigma\in{\bar {\frak H}}^*$, let ${\bar
h}_{\langle\sigma\rangle}\in{\bar
{\frak H}}$ be such that $\tau({\bar
h}_{\langle\sigma\rangle})={\bar
J}(\tau,\sigma)$
for all $\tau\in{\bar {\frak H}}^*$.
Then ${\bar h}_i={\bar \epsilon}_i{\bar h}_{{\langle{\bar
\al}_i\rangle}}$
and
${\bar t}_i={\bar \epsilon}_i{\bar h}_{{\langle{\bar
\gamma}_i\rangle}}$.

\begin{lemma}\label{lemma:affineformula}
Let ${\bar {\frak D}}$ be a handy datum. Then the
following
hold for
${\bar {\frak G}}({\bar A},{\bar
I}^{{\rm odd}})$.
\par
(1) $({\rm ad}{\bar E}_i)^{1-{\bar a}_{ij}}{\bar E}_j=0$
for $i\in{\bar I}^{{\rm pos}}$ and $j\in{\bar
I}\setminus\{i\}$. \par
(2) $[{\bar E}_i,{\bar E}_j]=0$ if ${\bar a}_{ij}=0$.
In particular it follows that if ${\bar a}_{ii}=0$,
then $[{\bar E}_i,{\bar E}_i]=0$ and  $({\rm ad}{\bar
E}_i)^2
X=0$ for any homogeneous element $X$ of ${\bar {\frak
G}}$.
\par
(3) $[{\bar E}_j,[[{\bar E}_i,{\bar E}_j],{\bar E}_m]]=0$
if
${\bar a}_{jj}=0$ and $-{\bar a}_{ji}={\bar a}_{jm}\ne 0$.
\par
(4) $[[{\bar E}_i,{\bar E}_j],{\bar E}_r]=[[{\bar
E}_i,{\bar E}_r],{\bar
E}_j]$
if ${\bar a}_{ii}=2$, ${\bar a}_{jj}={\bar a}_{rr}=0$,
${\bar a}_{ij}={\bar a}_{ji}={\bar a}_{ir}={\bar
a}_{ri}=-1$ and
${\bar a}_{jr}={\bar
a}_{rj}=2$.
\par (5) The same formulas as (1)-(4) with ${\bar F}_i$'s
in place of
${\bar E}_i$'s
hold. \par
(6) There exists a super-symmetric invariant form
${\bar J}:{\bar {\frak G}}\times{\bar {\frak G}}
\rightarrow\BBbb{C}$ such
that
${\bar J}({\bar h}_{{\langle\sigma\rangle}},{\bar
h}_{{\langle\tau\rangle}})={\bar
J}(\sigma,\tau)$ and
${\bar J}({\bar E}_i,{\bar F}_j)=\delta_{ij}{\bar
\epsilon}_i$.
(By abuse of notation, we use the same symbol ${\bar
J}$ for the bilinear forms on ${\bar {\frak G}}$ and
${\bar {\frak H}}^*$.)
\end{lemma}
This can be checked directly
(see also \cite[Proposition~6.7.1]{yamane1}).
\newline\par
Let ${\bar {\frak D}}$ be a handy datum. Let
$\BBbb{C}[t,t^{-1}]$
be the Laurent polynomial algebra. Let
\begin{equation}
{\frak L}({\bar {\frak D}}):={\bar {\frak G}}({\bar
A},{\bar
I}^{{\rm
odd}})\otimes\BBbb{C}[t,t^{-1}]\oplus\BBbb{C}v\oplus\BBbb{C}w.
\nonumber
\end{equation}
We view ${\frak L}({\bar {\frak D}})$ as a Lie
superalgebra in the
following way.
The parity of $X\otimes t^n$ is the same as the one of $X$
for a homogeneous element $X$ of ${\bar {\frak G}}({\bar
A},{\bar
I}^{{\rm odd}})$; the parities of $v$ and $w$ are $0$.
The Lie super bracket of ${\frak L}({\bar {\frak D}})$ is
given by
\begin{eqnarray*}
\lefteqn{[X\otimes t^m+a_1v+b_1w,Y\otimes t^n+a_2v+b_2w]}
\\
& &=[X,Y]\otimes t^{m+n}+m\delta_{m+n,0}{\bar J}(X,Y)v+
b_1nY\otimes t^n-b_2mX\otimes t^m
\end{eqnarray*} for homogeneous elements $X$, $Y$ of
${\bar
{\frak G}}({\bar
A},{\bar
I}^{{\rm odd}})$. We shall also denote by ${\bar J}$
the invariant form on ${\frak L}({\bar {\frak D}})$
defined by
\begin{equation}
{\bar J}(X\otimes t^m+a_1v+b_1w,Y\otimes t^n+a_2v+b_2w)
=\delta_{m+n,0}{\bar J}(X,Y)+a_1b_2+b_1a_2.
\nonumber
\end{equation}

\subsection{Unfolding}

Let ${\cal D}$ be a QEBS with $l\geq 2$.
We define a map $k^\vee:\Pi\rightarrow \{1,2,3,4\}$
by the following. \newline\par
({\rm KV}1) If $J(\beta^\vee,\al)=-1$ and
$g(\al)=\emptyset$
or $\BBbb{Z}$, then
$k^\vee(\al)={\frac {k(\beta)} {k(\al)}}k^\vee(\beta)$.
\par
({\rm KV}2) If $J(\al^\vee,\beta)=-2$, $k(\al)=2k(\beta)$
and
$g(\al)=2\BBbb{Z}$, then
$k^\vee(\al)=2k^\vee(\beta)$.
\par
({\rm KV}3) If $J(\al^\vee,\beta)=-2$, $k(\al)=k(\beta)$
and $g(\al)=2\BBbb{Z}+1$ or $2\BBbb{Z}$, then
$k^\vee(\al)=k^\vee(\beta)=2$. \par
({\rm KV}4) If $J(\al^\vee,\beta)=-2$ and
$g(\al)=4\BBbb{Z}$
or $4\BBbb{Z}+2$, then $k^\vee(\al)=3$ and
$k^\vee(\beta)=2$. \par
({\rm KV}5) $k^\vee(\al)=4k^\vee(\beta)$
if and only if $4k(\al)=k(\beta)$.
\newline\newline
For the above ${\cal D}$, we define a handy datum
${\bar {\frak D}}_{{\cal D}}=({\bar A}_{{\cal D}},
{\bar I}_{{\cal D}}, {\bar I}_{{\cal D}}^{{\rm odd}},
{\bar \Pi}_{{\cal D}},{\bar \Pi}_{{\cal D}}^\vee)$ in the
following way. Let
\begin{equation}
{\bar I}_{{\cal D}}
:=\{(\al,x)\in\Pi\times\{1,2,3,4\}|1\leq x\leq
k^\vee(\al)\}.
\nonumber
\end{equation}
We define a square
matrix ${\bar A}_{{\cal D}}
=({\bar a}_{(\al,x),(\beta,y)})_{(\al,x),(\beta,y)\in{\bar
I}_{{\cal D}}}$
in the following way. \par
({\rm AD}1) If $J(\al^\vee,\beta)=0$, then
${\bar a}_{(\al,x),(\beta,y)}=0$. \par
({\rm AD}2) Let $\al\in\Pi$. If $g(\al)=\emptyset$ or
$\BBbb{Z}$ then
${\bar a}_{(\al,x),(\al,y)}=2\delta_{xy}$. If
$g(\al)=2\BBbb{Z}+1$,
then ${\bar a}_{(\al,x),(\al,y)}=3\delta_{xy}-1$.
If $g(\al)=2\BBbb{Z}$,
then ${\bar a}_{(\al,x),(\al,y)}=2-2\delta_{xy}$.
If $g(\al)=4\BBbb{Z}$ or $4\BBbb{Z}+2$, then
${\bar a}_{(\al,1),(\al,1)}={\bar a}_{(\al,3),(\al,3)}$
$={\bar a}_{(\al,1),(\al,2)}={\bar
a}_{(\al,2),(\al,1)}=0$,
${\bar a}_{(\al,2),(\al,2)}={\bar
a}_{(\al,1),(\al,3)}={\bar
a}_{(\al,3),(\al,1)}=2$,
${\bar a}_{(\al,2),(\al,3)}=-1$ and ${\bar
a}_{(\al,3),(\al,2)}=-2$ \par
({\rm AD}3) Assume $J(\beta^\vee,\al)=-1$. Then
\begin{eqnarray}\nonumber
\lefteqn{({\bar a}_{(\al,x),(\beta,y)},{\bar
a}_{(\beta,y),(\al,x)})} \\
&=&
\left\{\begin{array}{ll}
(0,0) & \mbox{if $k^\vee(\al)=4$, $k^\vee(\beta)=2$ and
$x-y\notin
2\BBbb{Z}$,
} \\
(0,0) & \mbox{if $k^\vee(\beta)\leq k^\vee(\al)\leq {\frac
3
2}k^\vee(\beta)$ and
$x\ne y$, } \\
(-2,-1) & \mbox{if $2\leq k^\vee(\al)\leq 3$,
$k^\vee(\beta)=2$,
${\bar a}_{(\al,x),(\al,x)}=0$} \\
& \mbox{\,\,
and $x=y$,
} \\
(-1,-1) & \mbox{if $g(\al)=2\BBbb{Z}+1$ and $x=y$,} \\
({\frac {k(\al)} {k(\beta)}}
J(\al^\vee,\beta),-1) & \mbox{otherwise.}
\end{array}\right. \nonumber
\end{eqnarray}\par
Let ${\bar
E}_{-(\al,x)}:={\bar
F}_{(\al,x)}\in {\frak G}_{{\bar {\frak D}}_{{\cal D}}}$,
and notice that the ${\bar h}_{(\al,x)}
\in {\frak G}_{{\bar {\frak D}}_{{\cal D}}}$ is
not
necessarily ${\bar h}_{{\langle(\al,x)\rangle}}$.
\begin{lemma}\label{lemma:representation}
Let ${\cal D}$ be a QEBS with $l\geq 2$.
Then there exists a unique homomorphism $\pi_{{\cal
D}}:{\frak g}_{{\cal
D}}
\rightarrow{\frak L}({\bar {\frak D}}_{{\cal D}})$
satisfying the
following
properties:
\newline\par $({\rm PD}1)$
\begin{equation}\nonumber
\pi_{{\cal D}}(E_{\pm\al})=
\left\{\begin{array}{ll} \nonumber
\displaystyle{\sum_{x=1}^{k^\vee(\al)}{\bar
E}_{\pm(\al,x)}} &
\mbox{if $g(\al)=\emptyset$, $\BBbb{Z}$ or $2\BBbb{Z}$,}
\\
\sqrt{2}({\bar E}_{\pm(\al,1)}+{\bar E}_{\pm(\al,2)}) &
\mbox{if $g(\al)=2\BBbb{Z}+1$,} \\
\sqrt{2}{\bar E}_{\pm(\al,2)}\pm {\frac {1}
{\sqrt{2}}}[{\bar
E}_{\pm(\al,1)},{\bar
E}_{\pm(\al,3)}] & \mbox{if $g(\al)=4\BBbb{Z}+2$,} \\
{\bar E}_{\pm(\al,1)}\pm [{\bar E}_{\pm(\al,3)},{\bar
E}_{\pm(\al,2)}] & \mbox{if $g(\al)=4\BBbb{Z}$,}
\end{array}\right.
\end{equation}

\begin{eqnarray}\nonumber
\lefteqn{\pi_{{\cal D}}(E_{\pm\al^*})} \\
&=&\left\{\begin{array}{ll} \nonumber
\displaystyle{\sum_{x=1}^{k^\vee(\al)}
\zeta_\al^{\pm(2x-1-k^\vee(\al))}
{\bar E}_{\pm(\al,x)}\otimes t^{\pm k(\al)}} &
\mbox{if $g(\al)=\emptyset$ or $2\BBbb{Z}$,}  \\
\displaystyle{\pm {\frac 1 4}\sum_{x=1}^{k^\vee(\al)}
\zeta_\al^{\pm(2x-1-k^\vee(\al))}
[{\bar E}_{\pm(\al,x)},{\bar E}_{\pm(\al,x)}]\otimes
t^{\pm k(\al)}} &
\mbox{if $g(\al)=\BBbb{Z}$,}  \\
\sqrt{-1}[{\bar E}_{\pm(\al,1)},{\bar
E}_{\pm(\al,2)}]\otimes t^{\pm 1}  &
\mbox{if $g(\al)=2\BBbb{Z}+1$,} \\
({\bar E}_{\pm(\al,1)}+ \sqrt{-1}[{\bar
E}_{\pm(\al,3)},{\bar
E}_{\pm(\al,2)}])\otimes t^{\pm 1} & \mbox{if
$g(\al)=4\BBbb{Z}+2$,} \\
\sqrt{2}({\bar E}_{\pm(\al,2)}+{\frac {\sqrt{-1}} 2}
[{\bar
E}_{\pm(\al,1)},{\bar
E}_{\pm(\al,3)}])\otimes t^{\pm 1} & \mbox{if
$g(\al)=4\BBbb{Z}$,}
\end{array}\right.
\end{eqnarray} where $\zeta_\al:=
\exp({\frac {\pi\sqrt{-1}}
{k^\vee(\al)}})$.
\par

\par $({\rm PD}2)$ There exists a
$\kappa\in\BBbb{C}\setminus\{0\}$
such that ${\bar J}(\pi_{{\cal D}}(h_\mu),\pi_{{\cal
D}}(h_\nu))
=\kappa
      J(\mu,\nu)$ for $\mu$, $\nu\in{\cal E}$.
\par $({\rm PD}3)$ $\pi_{{\cal D}}(h_a)=\kappa v$,
$\pi_{{\cal D}}(h_{\Lambda_a})=w$
and $\pi_{{\cal D}}(h_{\Lambda_\delta})=
J(\Lambda_\delta,\al_0)\sum_{i=1}^{k^\vee(\al_0)}
{\bar t}_{(\al_0,i)}$.
In particular, $\pi_{{\cal D}}(h_\sigma)\ne 0$
for all $\sigma\in{\cal E}$.
\end{lemma}

This can be proved directly by using
Lemma~\ref{lemma:affineformula}.
\newline\par {\it Proof of
Lemma~\ref{lemma:Cartannondeg}.}
The lemma follows immediately from
Lemma~\ref{lemma:representation}. \hfill
$\BBox$
\newline\par
Keep the notation as in Lemma~\ref{lemma:representation}.
We shall also
denote by
$J$ the invariant form on ${\frak g}_{{\cal D}}$ defined
by ${\frac {1}
{\kappa}}{\bar
J}$. \par

\section{Invariant form and universal central extension
(UCE)}
\subsection{Invariant form and a universal property}
Let ${\cal D}$ be a QEBS with $l\geq 2$. Following the
notation in
\cite{saito1}, we say that ${\cal D}$ is $A_l^{(1,1)}$
if the $(l+1)\times(l+1)$-matrix
$A=(J(\al^\vee,\beta))_{\al,\beta\in\Pi}$
is
${\bf A}_l^{(1)}$ (see also Subsec.~1.1 for the name ${\bf
A}_l^{(1)}$).
Notice that if ${\cal D}$ is $A_l^{(1,1)}$, then
$p(\al)=0$ and $g(\al)=\emptyset$ for all $\al\in\Pi$.

\begin{theorem}\label{theorem:invarianttheorem}
Let ${\cal D}$ be a QEBS with $l\geq 2$.
Assume that ${\cal D}$ is
not
${\it A}_l^{(1,1)}$.
Let ${\frak g}^\prime_{{\cal D}}$ be a Lie superalgebra
satisfying the following conditions. \newline\par
$({\rm UI}1)$ ${\frak g}^\prime_{{\cal D}}$ includes
${\frak h}_{{\cal D}}$
as a sub Lie superalgebra.\par
$({\rm UI}2)$
\begin{equation}\nonumber
{\frak g}^\prime_{{\cal D}}={\frak h}_{{\cal D}}\bigoplus
\Bigl(\bigoplus_{\nu\in R(k,g)}{\frak g}^\prime_{{\cal
D},\nu}\Bigr)\bigoplus
\Bigl(\bigoplus_{(m,n)\in \BBbb{Z}^{2,\prime}}{\frak
g}^\prime_{{\cal
D},m\delta+na}\Bigr),
\end{equation} and $\dim{\frak g}^\prime_{{\cal D},\nu}=1$
for $\nu\in R(k,g)$,
where ${\frak g}^\prime_{{\cal D},\sigma}
:=[X\in {\frak g}^\prime_{{\cal
D}}|[h,X]=\sigma(h)X\,(h\in{\frak h}_{{\cal
D}})]$.
\par
$({\rm UI}3)$ The ${\frak g}^\prime_{{\cal D}}$ is
generated by
${\frak h}_{{\cal D}}$ and ${\frak g}^\prime_{{\cal
D},\nu}$
with $\nu\in R(k,g)$.
\par
$({\rm UI}4)$ There exists an invariant form $J^\prime$ on
${\frak g}^\prime_{{\cal
D}}$ such that $J^\prime(h_\sigma,h_\tau)=J(\sigma,\tau)$
($\sigma$,
$\tau\in{\cal
E}$) and such that $\ker J^\prime\subset
\bigoplus_{(m,n)\in
\BBbb{Z}^{2,\prime}}{\frak g}^\prime_{{\cal
D},m\delta+na}$.\newline\newline
Then there exists an epimorphism $\eta:{\frak g}_{{\cal
D}}
\rightarrow{\frak g}^\prime_{{\cal D}}$
such that $\eta(h_\sigma)=h_\sigma$
($\sigma\in{\cal E}$) and
$J^\prime\circ(\eta\times\eta)=J$.
\end{theorem}
 
{\it Proof.}
Using the same argument as in the proof of
\cite[Theorem~2.2]{kac1},
we can choose non-zero elements $E^\prime_\rho$ of ${\frak
g}^\prime_{{\cal
D},\rho}$
($\rho\in R(k,g)$)
so that $[E^\prime_\rho,E^\prime_{-\rho}]=h_{\rho^\vee}$.
By comparing (\ref{eqn:valueofn}) with
the equalities in ({\rm SR}6-7),
we can normalize the elements
$E^\prime_\mu$'s with $\mu\in{\cal B}$
so that they and the $h_\sigma$'s ($\sigma\in{\cal E}$)
satisfy the relations
({\rm SR}1-9).
Then the theorem follows from
Theorem~\ref{theorem:main}
and the existence of the $n_\mu$'s
with $\mu\in\Pi\cup\Pi^*$.
\hfill $\BBox$

\newtheorem{corollary}{Corollary}[section]
\begin{corollary}\label{corollary:galfourz}
Let ${\cal D}=\EPiakg$ be a QEBS with $l\geq 2$.
Assume $\al\in\Pi^B$ to be such that 
$g(\al)=4\BBbb{Z}$.
Let ${\cal D}^\prime={\cal D}({\cal E},\Pi,a,k,g^\prime)$
be the QEBS obtained from the ${\cal D}$ by replacing $g$
by $g^\prime$ such that $g^\prime(\al)=4\BBbb{Z}+2$ and
$g^\prime(\beta)=g(\beta)$
($\beta\ne\al$). Then there exists an isomorphism
$\xi:{\frak g}_{{\cal D}}\rightarrow
{\frak g}_{{\cal D}^\prime}$
such that $\xi (h_{\al})=h_{\al^*}$, $\xi
(h_{\beta})=h_{\beta}$
($\beta\ne\al$), $\xi (h_a)=h_a$,
$\xi (h_{\Lambda_\al})=h_{\Lambda_\al}$
and
$\xi (h_{\Lambda_a})=h_{\Lambda_a-\Lambda_\al}$,
where $\Lambda_\al\in\BBbb{C}\Pi\oplus\BBbb{C}a$
is such that $J(\Lambda_\al,\Lambda_\al)=0$ and
$J(\Lambda_\al,\gamma)=\delta_{\al\gamma}$
for $\gamma\in\Pi$.
\end{corollary}

This can be proved easily by using
Theorem~\ref{theorem:invarianttheorem}.

\subsection{A Lie algebra with the quantum tori}
Here we recall a Lie algebra studied in
\cite{BermanGaoKrylyuk}.
Let $q\in\BBbb{C}\setminus\{0\}$.
Let $\BBbb{C}_q=\BBbb{C}_q[s^{\pm 1},t^{\pm 1}]$ be the
$\BBbb{C}$-algebra
defined by generators $s^{\pm 1}$, $t^{\pm 1}$ and
defining relations
$ts=qst$. Let ${\rm M}_{l+1}(\BBbb{C}_q)$ be the
$\BBbb{C}_q$-algebra
of the $(l+1)\times (l+1)$-matrices over $\BBbb{C}_q$.
Let  ${\widehat {\rm M}}_{l+1}(\BBbb{C}_q):={\rm
M}_{l+1}(\BBbb{C}_q)
\oplus\BBbb{C}c_1\oplus\BBbb{C}c_2\oplus\BBbb{C}d_1\oplus\BBbb{C}d_2$.
We regard ${\widehat {\rm M}}_{l+1}(\BBbb{C}_q)$ as a
$\BBbb{C}$-Lie
algebra
by
\begin{eqnarray*}\nonumber
\lefteqn{[s^{x_1}t^{x_2}E_{ij},s^{y_1}t^{y_2}E_{mn}]} \\
&=& s^{x_1+y_1}t^{x_2+y_2}(\delta_{jm}q^{x_2y_1}E_{in}
-\delta_{in}q^{x_1y_2}E_{mj}) \\
&
&\,\,+\delta_{x_1+y_1,0}\delta_{x_2+y_2,0}q^{x_2y_1}(x_1c_1+x_2c_2),
\end{eqnarray*}
\begin{equation}\nonumber
[c_i,s^{x_1}t^{x_2}E_{mn}]=0,\,
[d_i,s^{x_1}t^{x_2}E_{mn}]=x_is^{x_1}t^{x_2}E_{mn}\,
\end{equation}
and
\begin{equation}\nonumber
[c_i,c_j]=[c_i,d_j]=[d_i,d_j]=0.
\end{equation}
We define a symmetric invariant form ${\bar J}_q$ on
${\widehat {\rm
M}}_{l+1}(\BBbb{C}_q)$ by
\begin{equation}\nonumber
\left\{\begin{array}{l}
{\bar
J}_q(s^{x_1}t^{x_2}E_{ij},s^{y_1}t^{y_2}E_{mn})
=\delta_{x_1+y_1,0}\delta_{x_2+y_2,0}q^{x_2y_1}, \\
{\bar J}_q(s^{x_1}t^{x_2}E_{ij},c_i)={\bar
J}_q(s^{x_1}t^{x_2}E_{ij},d_i)=0, \\
{\bar J}_q(c_i,c_j)={\bar J}_q(d_i,d_j)=0,\,{\bar
J}_q(c_i,d_j)=\delta_{ij}.
\end{array}\right.
\end{equation}
     \par
Let ${\cal D}$ be ${\it
A}_l^{(1,1)}$.
Let ${\frak g}_{{\cal D}}^q$
be the Lie algebra defined be the generators
(\ref{eqn:generators})
and the defining relations
obtained from the ones in ({\rm SR}1-9) by replacing
$[E_{\pm\al^*_0},E_{\pm\al_l}]=[E_{\pm\al_0},E_{\pm\al^*_l}]$
of ({\rm SR}6-7)
with
\begin{equation}\nonumber
q^{\pm
1}[E_{\pm\al^*_0},E_{\pm\al_l}]=[E_{\pm\al_0},E_{\pm\al^*_l}].
\end{equation} Then there exists a unique homomorphism
$\pi_{{\cal D}}^q:{\frak g}_{{\cal
D}}^q\rightarrow{\widehat {\rm
M}}_{l+1}(\BBbb{C}_q)$
such that $\pi_{{\cal D}}^q(E_{\al_i})=E_{ii+1}$,
$\pi_{{\cal
D}}^q(E_{\al^*_i})=tE_{ii+1}$,
      $\pi_{{\cal D}}^q(E_{-\al_i})=E_{i+1i}$, $\pi_{{\cal
D}}^q(E_{-\al^*_i})=t^{-1}E_{i+1i}$ ($1\leq i\leq l$),
$\pi_{{\cal D}}^q(E_{\al_0})=sE_{l+11}$,
$\pi_{{\cal D}}^q(E_{\al^*_0})=stE_{l+11}$,
$\pi_{{\cal D}}^q(E_{-\al_0})=s^{-1}E_{1l+1}$,
$\pi_{{\cal D}}^q(E_{-\al^*_0})=qs^{-1}t^{-1}E_{1l+1}$,
$\pi_{{\cal D}}^q(h_{\Lambda_\delta})=d_1$ and $\pi_{{\cal
D}}^q(h_{\Lambda_a})=d_2$. We see that similar results to
Theorems
\ref{theorem:main} and \ref{theorem:invarianttheorem}
also hold for ${\frak g}_{{\cal D}}^q$ with $\pi_{{\cal
D}}^q$. Especially we have the following.

\begin{theorem}\label{theorem:qinvarianttheorem}
Let ${\cal D}$ be a QEBS with $l\geq 2$.
Assume that ${\cal D}={\it A}_l^{(1,1)}$.
Let ${\frak g}^\prime_{{\cal D}}$ be a Lie superalgebra
satisfying the same conditions as the $\mbox{({\rm
UI}1-4)}$ 
in Theorem~\ref{theorem:invarianttheorem}. 
Then there exist a $q\in\BBbb{C}\setminus\{0\}$
and
an epimorphism $\eta:{\frak g}_{{\cal
D}}^q
\rightarrow{\frak g}^\prime_{{\cal D}}$
such that $\eta(h_\sigma)=h_\sigma$
($\sigma\in{\cal E}$) and
$J^\prime\circ(\eta\times\eta)=J$.
\end{theorem}

Using \cite[Chap.~I Theorem
1.29(d)]{Allisonetal},
by
Theorems~\ref{theorem:invarianttheorem}
and \ref{theorem:qinvarianttheorem}, we have
the following.
\begin{theorem}\label{theorem:EALA}
Let ${\cal D}=\EPiakg$ be a QEBS with $l\geq 2$.
We assume that $g(\al)=\emptyset$ or
$2\BBbb{Z}+1$ for every $\al\in\Pi$.
We also keep the notation in \cite[Chap.~I
\S1]{Allisonetal}.
Let $({\cal L},(\,,\,),{\cal H})$ be an extended affine
Lie algebra (EALA)
in the sense of \cite[Chap.~I Definition
1.33]{Allisonetal}. Assume that
there exists an isometric monomorphism
$\varphi:{\cal E}\rightarrow{\cal H}^*$ such
that the set $\varphi (\Rkg)$ coincides with the set
of the non-isotropic roots of the ${\cal L}$.
Then there exists a homomorphism
${\cal F}:{\bar {\frak g}}_{\cal D}\rightarrow{\cal L}$
such that ${\cal F}(h_\sigma)=t_{\varphi(\sigma)}$
($\sigma\in{\cal E}$)
and ${\cal F}({\frak g}_\mu)={\cal L}_{\varphi(\mu)}$
($\mu\in\Rkg$), where we let ${\bar {\frak g}}_{\cal
D}:={\frak
g}_{\cal D}^q$
for some $q\in\BBbb{C}\setminus\{0\}$
if ${\cal D}=A_l^{(1,1)}$, and, otherwise,
we let ${\bar {\frak g}}_{\cal D}:={\frak g}_{\cal D}$.
(See \cite[Chap.~I (1.2) and (1.8)]{Allisonetal}
for the symbols ${\cal L}_{\varphi(\mu)}$ and
$t_{\varphi(\sigma)}$. )
\end{theorem}

\subsection{Universal central extension (UCE)}
We first recall the definition of the universal central
extension (UCE for short)
of a Lie superalgebra. See \cite{IoharaKoga} for more
detail.
Let ${\frak a}={\frak a}_0\oplus{\frak a}_1$ be a Lie
superalgebra.
Let ${\Der}({\frak a}):=[{\frak a},{\frak a}]$.
We say that ${\frak a}$ is perfect if ${\Der}({\frak
a})={\frak a}$.
We say that a Lie superalgebra epimorphism
$P:{\frak u}={\frak u}_0\oplus{\frak u}_1
\rightarrow {\frak a}$ is a central extension
if $[\ker P,{\frak u}]=\{0\}$ and
$\ker P=({\frak u}_0\cap\ker P)\oplus({\frak u}_1\cap\ker
P)$ (we do not assume $\ker P\subset
{\frak u}_0$). We say that a central extension
$V:{\frak u}\rightarrow{\frak a}$
is a UCE if ${\frak u}={\Der}({\frak u})$ and if
for any central
extension
$W:{\frak b}\rightarrow{\frak a}$, there exists a
homomorphism
$M:{\frak u}\rightarrow{\frak b}$ such that $W\circ M=V$.
Notice that if $P:{\frak u}
\rightarrow {\frak a}$ is a central extension and if $x$,
$y\in{\frak a}$
are homogeneous elements,
then there exists a unique $z\in{\frak u}$
such that $z\in[P^{-1}(\{x\}),P^{-1}(\{y\})]$; we denote
the $z$
by $N(P,x,y)$. Notice that if $x\in{\frak a}_i$ and
$y\in{\frak a}_i$,
then $N(P,x,y)\in{\frak u}_{i+j}$.
\par
The following lemma seems to be trivial, but
we give it to use it in the proof of
Theorem~\ref{theorem:centralextension}.
\begin{lemma}\label{lemma:easylemma}
Let ${\frak a}={\frak a}_0\oplus{\frak a}_1$ be a Lie
superalgebra
such that ${\frak a}={\rm
Der}({\frak
a})\oplus\BBbb{C}k_1\oplus\cdots\oplus\BBbb{C}k_n$,
$k_i\in{\frak a}_0$, $[k_i.k_j]=0$, and
${\Der}({\frak a})=\oplus_{x\in\BBbb{C}^n}{\frak
a}^\prime_x$,
where ${\frak a}^\prime_x=\{X\in{\Der}({\frak
a})|[k_i,X]=x_iX\}$
(here $x=(x_1,\ldots,x_n$)).
Assume that ${\frak a}$ is presented by generators
$k_i$ ($1\leq i\leq n$)
$a_p\in{\Der}({\frak a})$ ($p\in P$) with
$a_p\in{\frak a}^\prime_{x(p)}$
for some
$x(p)=(x(p)_1,\ldots,x(p)_n)\in\BBbb{C}^n$ and defining
relations $f_t=0$ ($t\in T$) and $[k_i,a_p]=x(p)_ia_p$,
$[k_i.k_j]=0$,
where $f_t$'s are assumed to be expressed only by the
elements $a_p$ ($p\in
P$) and to be homogeneous with respect to the $a_p$'s.
Then ${\Der}({\frak a})$ is also presented by the
generators
$a_p$ ($p\in P$) and the defining relations $f_t=0$ ($t\in
T$).
\end{lemma} \par
{\it Proof.} Let ${\frak
c}:=\BBbb{C}k_1\oplus\cdots\oplus\BBbb{C}k_n$.
Let ${\frak b}$ be the Lie superalgebra generated by the
generators
$b_p$ ($p\in P$) and the defining relations $g_t=0$ ($t\in
T$),
where $g_t$ is expressed by replacing  $a_p$ of $f_t$ with
$b_p$.
For $x\in\BBbb{C}^n$, let ${\frak b}_x$ be the subspace of
${\frak b}$
spanned by the elements
$({\rm ad} b_{p_1})\cdots({\rm ad} b_{p_{r-1}})b_{p_r}$
with $x(p_1)+\cdots +x(p_r)=x$.
Then ${\frak b}=\oplus_{x\in\BBbb{C}^n}{\frak b}_x$.
We can define a Lie superalgebra ${\frak d}={\frak
b}\oplus{\frak c}$
by $[b+\sum y_ik_i,b^\prime+\sum y^\prime_ik_i]
=[b,b^\prime]+(\sum y_ix^\prime_i)b^\prime-(\sum
y^\prime_ix_i)b$,
where $b\in{\frak b}_x$ and $b^\prime\in{\frak
b}_{x^\prime}$.
We see that there exists an isomorphism
$\Phi:{\frak d}\rightarrow{\frak a}$
such that $\Phi(b_p)=a_p$ and $\Phi(k_i)=k_i$. \hfill
$\BBox$ \newline\par
Now we give a UCE theorem.

\begin{theorem}\label{theorem:centralextension}
Let ${\cal D}$ be a QEBS with $l\geq 2$.
Recall the ${\bar {\frak g}}_{{\cal D}}$
from Theorem~\ref{theorem:EALA}.
Let $\varrho:{\Der}({\bar {\frak g}}_{{\cal D}})
\rightarrow{\frak a}$
be a central extension.
Then the $\varrho$
is a UCE.
\end{theorem} \par
{\it Proof.} We first assume that
${\bar {\frak g}}_{{\cal D}}={\frak g}_{{\cal D}}$.
Let $f:{\frak b}\rightarrow{\frak a}$ be a central
extension.
For $\mu\in{\cal B}_+$, let
\begin{equation}\nonumber
h^\prime_{\mu^\vee}:=N(f,\varrho(E_\mu),\varrho(E_{-\mu}))
\quad\mbox{and}\quad
E^\prime_{\pm\mu}:=N(f,\varrho(\pm{\frac {1} {2}}
h_{\mu^\vee}),\varrho(E_{\pm\mu})).
\end{equation} By Lemma~\ref{lemma:easylemma}, it suffices
to show
that these elements of the ${\frak b}$ satisfy the
equalities
in the ({\rm SR}1-9).
\par
We first show the ({\rm SR}4). Let $\mu\in{\cal B}_+$.
Notice that
\begin{equation}\nonumber
\{f(E^\prime_{\pm\mu})\}=
\{f([f^{-1}(\{\varrho(\pm{\frac {1} {2}}
h_{\mu^\vee})\}) 
,f^{-1}(\{\varrho(E_{\pm\mu})\})])\}
=\{\varrho(E_{\pm\mu})\}.
\end{equation} 
Hence we have
\begin{equation}\label{eqn;fepmm}
 f(E^\prime_{\pm\mu})=\varrho(E_{\pm\mu})\quad
\mbox{and}\quad 
E^\prime_{\pm\mu}\in f^{-1}(\{\varrho(E_{\pm\mu})\}),
\end{equation}
which implies the ({\rm SR}4), as desired.
\par We show the ({\rm SR}2-3). For $\mu\in{\cal B}_+$, 
by (\ref{eqn;fepmm}) and the ({\rm SR}4),
we have $f(h^\prime_{\mu^\vee})
=[f(E^\prime_\mu),f(E^\prime_{-\mu})]=\varrho(h_{\mu^\vee})$,
and, by (\ref{eqn;fepmm}), we have $E^\prime_{\pm\mu}=
[\pm{\frac {1} {2}}
h^\prime_{\mu^\vee},E^\prime_{\pm\mu}]$.
Hence, for $\mu$, $\nu\in{\cal B}_+$, we have
\begin{equation}\nonumber
[h^\prime_{\nu^\vee},E^\prime_{\pm\mu}]
=[h^\prime_{\nu^\vee},[\pm{\frac {1} {2}}
h^\prime_{\mu^\vee},E^\prime_{\pm\mu}]]
=\pm{\frac {1} {2}}[h^\prime_{\mu^\vee},[
h^\prime_{\nu^\vee},E^\prime_{\pm\mu}]]
=J(\nu^\vee,\pm\mu)E^\prime_{\pm\mu}, 
\end{equation} and
$
[h^\prime_{\nu^\vee},h^\prime_{\mu^\vee}]=
[h^\prime_{\nu^\vee},[E^\prime_\mu,E^\prime_{-\mu}]]=0 
$,
as desired. \par
For $\mu\in{\cal B}_+$, let
$h^\prime_{(-\mu)^\vee}:=-h^\prime_{\mu^\vee}$.
\par
We show the ({\rm SR}5-9). For 
$(\mu,\nu)\in ({\cal B}\times{\cal B})^\prime$
and $y$ $z\in\BBbb{C}$, by ({\rm SR}3), we have
\begin{equation}\nonumber
0=[yh^\prime_{\mu^\vee}+zh^\prime_{\nu^\vee},
({\rm ad}E^\prime_\mu)^{x_{\mu,\nu}}E^\prime_\nu]
=J(y\mu^\vee+z\nu^\vee,
x_{\mu,\nu}\mu+\nu)({\rm
ad}E^\prime_\mu)^{x_{\mu,\nu}}E^\prime_\nu,
\end{equation} which implies 
$({\rm ad}E^\prime_\mu)^{x_{\mu,\nu}}E^\prime_\nu=0$.
Hence we have the ({\rm SR}5). Similarly we have ({\rm
SR}6-9).
\par
Finally we show the ({\rm SR}1). For
$\al\in\Pi$, let 
\begin{equation}\nonumber
h^\prime_{a,\al}:={\frac {J(\al,\al)c(\al)} {2k(\al)}}
(c(\al)h^\prime_{(\al^*)^\vee}-h^\prime_{\al^\vee}). 
\end{equation}
It suffices to show that
$h^\prime_{a,\al}=h^\prime_{a,\beta}$
for all $\al$, $\beta\in\Pi$.
By ({\rm SR}3-5), we see that for each $\mu\in{\cal B}$, 
$E^\prime_\mu$ is locally nilpotent,
so $n_\mu=n_{E^\prime_\mu}\in{\rm Aut}({\frak b})$
can be defined in the same way as in (2.3). Let
$(\al,\beta,y)\in{\cal A}$. By the ({\rm SR}6-7), we have
$n_{\al^*}E^\prime_{\pm\beta}=n_\al
E^\prime_{\pm\beta^*}$.
Hence $n_{\al^*}h^\prime_{\beta^\vee}=n_\al
h^\prime_{(\beta^*)^\vee}$.
Since $n_{\al^*}h^\prime_{\beta^\vee}
=h^\prime_{\beta^\vee}-J(\al^*,\beta^\vee)
=h^\prime_{\beta^\vee}-c(\al)J(\al,\beta^\vee)h^\prime_{(\al^*)^\vee}$
and $n_\al
h^\prime_{(\beta^*)^\vee}=h^\prime_{(\beta^*)^\vee}
-J(\al,\beta^\vee)h^\prime_{\al^\vee}$, we have 
$h^\prime_{a,\al}=h^\prime_{a,\beta}$, as desired.
 \par
We see that the case where ${{\cal D}}=A_l^{{(1,1)}}$ and
${\bar {\frak g}}_{{\cal D}}={\frak g}_{{\cal D}}^q$ with
$q\ne 1$
can also be treated similarly.
\hfill $\BBox$ \newline\par
Keep the notation as above.
Denote by ${\bar \pi}_{{\cal D}}$
the homomorphism $\pi_{{\cal D}}:{\frak g}_{{\cal D}}
\rightarrow{\frak L}({\bar {\frak D}}_{{\cal D}})$ in
Lemma~\ref{lemma:representation}
if ${\frak g}_{{\cal D}}={\bar {\frak g}}_{{\cal D}}$;
and, otherwise, let ${\bar \pi}_{{\cal D}}$
denote the homomorphism $\pi_{{\cal D}}^q:{\frak g}_{{\cal
D}}^q
\rightarrow {\hat M}_{l+1}(\BBbb{C}_q)$.
Let $\varpi:{\rm Im}{\bar \pi}_{{\cal D}}\rightarrow
{\rm Im}{\bar \pi}_{{\cal D}}/\pi_{{\cal
D}}(\BBbb{C}h_\delta
\oplus\BBbb{C}h_a)$ be the
natural projective map. Define
the epimorphism ${\hat \pi}_{{\cal D}}:{\Der}({\bar
{\frak g}}_{{\cal D}})
\rightarrow(\varpi\circ{\bar \pi}_{{\cal D}})({\rm
Der}({\bar {\frak
g}}_{{\cal D}}))$ by
${\hat \pi}_{{\cal D}}=(\varpi\circ{\bar \pi}_{{\cal
D}})_{|{\rm
Der}({\bar {\frak g}}_{{\cal
D}})}$.
By Theorem~\ref{theorem:centralextension}, we see that the
${\hat
\pi}_{{\cal D}}$
is a UCE.
In particular, we have the following.

\begin{corollary}\label{corollary:MEY}
(1) The ${\hat \pi}_{{\cal D}}$ is a UCE. \par
(2) Keep the notation as in Theorem~\ref{theorem:EALA}.
Then the ${\cal F}_{|{\Der}({\bar {\frak g}}_{{\cal
D}})}$
$:{\Der}({\bar {\frak g}}_{{\cal D}})\rightarrow$
${\cal F}({\Der}({\bar {\frak g}}_{{\cal D}}))$
is a UCE. (Notice that the ${\cal F}({\Der}
({\bar {\frak g}}_{{\cal D}}))$ is the core of
the EALA ${\cal L}$ in the sense of
\cite[Chap.~I Definition 2.20]{Allisonetal}.)
\end{corollary}
{\it Proof.} The statement (1) follows from the argument
in the paragraph above the corollary.
The statement (2) also follows from
Theorem~\ref{theorem:centralextension}.
\hfill $\BBox$ \newline\par
By Corollary~\ref{corollary:MEY}, we see the following.
If $g(\al)=\emptyset$ for all $\al\in\Pi$ and if the
Cartan matrix $A$ is ${\bf X}_l^{(1)}$
for some ${\bf X}={\bf A},\ldots,{\bf G}$, then ${\rm
Der}({\frak g}_{{\cal
D}})$
is isomorphic to the (2-variable) toroidal Lie algebra in
the sense of
\cite{MoodyRaoYokonuma}. If the $A={\bf A}_l^{(1)}$, i.e.,
${\cal D}=A_l^{(1,1)}$, then ${\rm
Der}({\frak g}_{{\cal
D}}^q)$ is isomorphic to the Steinberg Lie algebra
$st_{l+1}(\BBbb{C}_q)$ (see \cite{BermanGaoKrylyuk}
for the term and symbol).

\section{Elliptic root base}

\subsection{Weakening the condition of ({\rm SR}5)}
Here we show that some of the ({\rm SR}5) are not
necessary. Let ${\cal D}=\EPiakg$ be a QEBS with $l\geq 2$. 
Let
$(\Pi\times\Pi)^\prime:=\{(\al,\beta)\in\Pi\times\Pi\,|\,
\al\ne\beta\}$.
Define a subset $(\Pi\times\Pi)^\sharp$ of
$(\Pi\times\Pi)^\prime$ by
\begin{eqnarray*}
(\Pi\times\Pi)^\sharp&:=&\{(\al,\beta)\in
(\Pi\times\Pi)^\prime\,|\,
k(\al)=k(\beta),\,J((\al^*)^\vee,\beta)=-1\} \\
& &\cup \{(\al,\beta)\in (\Pi\times\Pi)^\prime\,|\,
k(\al)<k(\beta),\,J(\beta^\vee,\al)=-1\}.
\end{eqnarray*}
For a subset ${\cal X}$ of  $({\cal B}\times{\cal
B})^\prime$,
let $PM({\cal
X}):=\{(\varepsilon_1\mu,\varepsilon_2\nu)\,|\,
(\mu,\nu)\in{\cal X},\,
\varepsilon_1,\,\varepsilon_2\in\{1,-1\}\}$.
Define a subset $({\cal B}\times{\cal B})^\sharp$
of $({\cal B}\times{\cal B})^\prime$ by
\begin{eqnarray*}
\lefteqn{({\cal B}\times{\cal B})^\sharp} \\
&:=&\{(\mu,\nu)\in({\cal B}\times{\cal B})^\prime\,|\,
J(\mu,\nu)=0 \} \\
& &\,\cup PM(\cup_{(\al,\beta)\in(\Pi\times\Pi)^\sharp}
\{(\al,\beta),(\beta,\al),
(\al^*,\beta),(\beta,\al^*),
(\al,\beta^*)
\}
) \\
& &\,\cup PM(\cup_{\al\in\Pi}
\{(\al,\al^*),
(\al^*,\al)
\}
).
\end{eqnarray*}

Then we consider the following.
\newline
\par
$({\rm SR}5^\prime)$\quad
$({\rm ad}E_\mu)^{x_{\mu,\nu}}E_\nu=0$ \quad
if $(\mu,\nu)\in( {\cal B}\times{\cal B})^\sharp$.
\begin{theorem}\label{theorem:reduced}
Let ${\cal D}$ be a QEBS with $l\geq 2$. Let
${\frak g}_{{\cal D}}^\sharp$ be the Lie superalgebra
defined by the same generators as in
(\ref{eqn:generators})
with the same parities as in
(\ref{eqn:generatorsparities})
and by the defining relations $\mbox{{\rm ({\rm
SR}1-4)}}$,
$\mbox{{\rm ({\rm SR}$5^\prime$)}}$, $\mbox{{\rm ({\rm
SR}6-9)}}$.
Define the homomorphism $\Theta:{\frak g}_{{\cal
D}}^\sharp
\rightarrow{\frak g}_{{\cal D}}$
by $\Theta(h_\sigma)=h_\sigma$ $(\sigma\in{\cal E})$
and $\Theta(E_\mu)=E_\mu$ $(\mu\in{\cal B})$.
Then $\Theta$ is an isomorphism.
\end{theorem} \par
{\it Proof.} Let $(\al,\beta)\in(\Pi\times\Pi)^\sharp$ and
let $T:=\{\al,\al^*,\beta\}$. Let
$({\frak g}_{{\cal D}}^\sharp)^{(T)}$ be
the sub Lie superalgebra of ${\frak g}_{{\cal D}}^\sharp$
generated by the $h_\sigma$ ($\sigma\in{\cal E}$)
and $E_\mu$ ($\mu\in T\cup-T$). Then we see that
for each $\mu\in T\cup -T$, the $E_\mu$ is a locally
nilpotent element
of $({\frak g}_{{\cal D}}^\sharp)^{(T)}$, so
$n_\mu\in{\rm Aut}(({\frak g}_{{\cal D}}^\sharp)^{(T)})$
can
be defined
in the same way as in (\ref{eqn:nmu}).
By ({\rm SR}6-7), we have
$E_{\pm\beta^*}=n_\al^{-1}n_\al E_{\pm\beta^*}
=n_\al^{-1}n_{\al^*}E_{\pm\beta}\in({\frak g}_{{\cal
D}}^\sharp)^{(T)}$.
It follows from
$({\rm SR}5^\prime)$ and ({\rm SR}6-7) that
\begin{equation}\label{eqn;alalstarbeta}
[E_{\pm\al},[E_{\pm\al^*},E_{\pm\beta}]]
=({\rm ad}E_{\pm\al})^{c(\al)+1}E_{\pm\beta^*}=0
\quad\mbox{if $k(\al)=k(\beta)$.}
\end{equation}
By $({\rm SR}5^\prime)$, we have
\begin{equation}\label{eqn:mpalpmbetastar}
[E_{\mp\al},E_{\pm\beta^*}]=0.
\end{equation}
It follows from (\ref{eqn;alalstarbeta}) and 
({\rm SR}1-4,$5^\prime$,8-9), that
\begin{eqnarray}
\lefteqn{[E_{\mp\al^*},E_{\pm\beta^*}]}\label{eqn:alstarbetastar}
\\
&
\sim &
n_\al^{-1}[({\rm ad}
E_{\pm\al})^{J(\al^\vee,\al^*)}E_{\mp\al^*},
({\rm ad}
E_{\pm\al^*})^{-J((\al^*)^\vee,\beta)}E_{\pm\beta}]
\nonumber \\
&
\sim & n_\al^{-1}({\rm ad} E_{\pm\al})^{J(\al^\vee,\al^*)}
({\rm ad}
E_{\pm\al^*})^{-J((\al^*)^\vee,\beta)-1}E_{\pm\beta}=0.
\nonumber
\end{eqnarray}
Since the $E_{\pm\al^*}$ and $E_{\pm\beta^*}$
are locally nilpotent in $({\frak g}_{{\cal
D}}^\sharp)^{(T)}$,
by (\ref{eqn:mpalpmbetastar}) and
(\ref{eqn:alstarbetastar}),
we see that
$({\rm ad}E_\mu)^{x_{\mu,\nu}}E_\nu=0$ 
for $(\mu,\nu)\in PM(\{(\beta^*,\al),(\al^*,\beta^*),
(\beta^*,\al^*)\})$.
This completes the proof. \hfill $\BBox$.
\newline\par
We notice that there are also redundant equalities
in ({\rm SR}$5^\prime$), which can be seen by the
calculations
in \cite[\S2.3]{yamane2} and in Proof of
Lemma~\ref{lemma:embedinglemma}.
We also notice that the condition of
the equality $[E_{\al^*},E_{-\beta^*}]=0$
of \cite[({\rm S}4)]{yamane2}
can be weakened to be the one that $\al$,
$\beta\in\Pi_{{\rm af}}$
with $I(\al,\beta)=0$,
which can be seen the same calculations
as in (\ref{eqn:alstarbetastar}).

\subsection{Relations for the elliptic root basis}

Let ${\cal D}=\EPiakg$ be a QEBS with $l\geq 2$.
Here by extending the notion given in [S], we introduce an
{\it
elliptic root basis} of $\Rkg$ for the ${\cal D}$.
Recall the $\delta\in\BBbb{Z}_+\Pi$ from Subsec.~2.2.
For $\al\in\Pi$, let $x_\al\in\BBbb{Z}_+$ ($\al\in\Pi$)
be the coefficient of $\delta$, i.e.,
$\delta=\sum_{\al\in\Pi}x_\al\al$.
Let $m_\al:={\frac {c(\al)I(\al,\al)x_\al}
{k(\al)}}$.
Let $m_{{\rm max}}:={\rm max}\{m_\al |\al\in\Pi\}$ and
$\Pi_{{\rm max}}:=\{\al\in\Pi | m_\al =m_{{\rm max}}\}$.
For a subset $S$ of $\Pi$, let $S^*:=\{\al^*|\al\in S\}$.
Let $\Gamma(R,G):=\Pi\cup\Pi_{{\rm max}}^*$
(where $R$ and $G$ denote $\Rkg$ and $\BBbb{C}a$
respectively).
We call the $\Gamma(R,G)$ the {\it elliptic root basis} of
$\Rkg$.
For a subset $S$ of $\Pi$, let
$\Gamma(R,G;S):=\Gamma(R,G)\cap(S\cup
S^*)$.
Recall the Lie superalgebra ${\frak g}_{{\cal D}}$ from
Subsec.~2.1.
\begin{theorem}\label{theorem:ellipticgenerators}
Let ${\cal D}$ be a QEBS with $l\geq 2$.
Let ${\frak g}_{{\cal D}}^{\Gamma(R,G)}$ be
the Lie superalgebra defined by the generators
\begin{equation}\label{eqn:ellipticgenerators}
h_\sigma\,(\sigma\in{\cal E}),
E_\mu,\,E_{-\mu}\,(\mu\in\Gamma(R,G))
\end{equation} with the same parities as
(\ref{eqn:generatorsparities})
and the following defining relations.
\newline\par
$({\rm TSR}i)$ $(1\leq i\leq 4,\,i=5^\prime,\,
6\leq i\leq 9)$
The same relations as the ones among those in $({\rm
SR}i)$
expressed only by the same symbols as in
(\ref{eqn:ellipticgenerators}), \par
$({\rm TSR}10)$ \quad
$[E_{\pm\al^*},[E_{\pm\al},E_{\pm\beta}]]=0$,\quad
if $\Gamma(R,G;\{\al,\,\beta\})=\{\al,\al^*,\beta\}$ and
${\frac {J(\al^\vee,\beta)}
{J(\beta^\vee,\al)}}=c(\al)$,\par
$({\rm TSR}11)$ \quad
$[E_{\pm\al^*},[E_{\pm\al},E_{\pm\beta}]]=0$ and
$[[E_{\pm\al^*},E_{\pm\beta}],[E_{\pm\al},E_{\pm\beta}]]=0$,
\quad if
$\Gamma(R,G;\{\al,\,\beta\})=\{\al,\al^*,\beta\}$,
$g(\al)=\emptyset$
and
${\frac {J(\beta^\vee,\al)} {J(\al^\vee,\beta)}}=2$ or
$3$,
\par
$({\rm TSR}12)$ \quad
$[({\rm ad}E_{\pm\beta})^{-J(\beta^\vee,\al)}
E_{\pm\al},({\rm
ad}E_{\pm\beta^*})^{-J((\beta^*)^\vee,\gamma)}
E_{\pm\gamma}]=0$, \newline
if
$\Gamma(R,G;\{\al,\,\beta,\,\gamma\})=\{\al,\beta,\beta^*,\gamma\}$,
\newline\newline
where we assume the $\al$, $\beta$, $\gamma$ are distinct.
Then
there exists a unique isomorphism $\Xi:{\frak g}_{{\cal
D}}^{\Gamma(R,G)}
\rightarrow{\frak g}_{{\cal D}}^\sharp$ such that
$\Xi(h_\sigma)=h_\sigma$ and $\Xi(E_\mu)=E_\mu$
($\mu\in{\cal B}^{\Gamma(R,G)}$), where we set
${\cal B}^{\Gamma(R,G)}:=\Gamma(R,G)\cup -\Gamma(R,G)$.
\end{theorem} \par
{\it Proof.} If the ${\cal D}$ is an SQEBS, then the
theorem
can be proved by the same argument as that for
\cite[Theorem~4.1]{yamane2}
(See also the last paragraph in Subsec.~6.1).
\par Assume the ${\cal D}$ not to be an SQEBS.
Then, for each $\beta\in\Pi\setminus\Pi_{{\rm {max}}}$,
there exists a unique $\al_\beta\in\Pi\setminus\{\beta\}$
such that
$J(\beta,\al_\beta)\ne 0$, and we see that
$\al_\beta\in\Pi_{{\rm {max}}}$ and
$J(\beta^\vee,\al_\beta)J(\beta,\al_\beta^\vee)=1$ or
$2$.
By Theorems \ref{theorem:main} and \ref{theorem:reduced},
we see that the homomorphism
$\Xi$ in
the statement exists.
By the same argument as in the proof of
Theorem~\ref{theorem:reduced}, we see that 
for every $(\mu,\nu)
\in(({\cal B}\times{\cal B})^\prime\cap({\cal
B}^{\Gamma(R,G)}
\times{\cal B}^{\Gamma(R,G)}))$,
the equality $({\rm ad}E_\mu)^{x_{\mu,\nu}}E_\nu=0$
holds in 
${\frak g}_{{\cal D}}^{\Gamma(R,G)}$.
Hence the elements $E_\mu$
($\mu\in{\cal B}^{\Gamma(R,G)}$) of ${\frak g}_{{\cal
D}}^{\Gamma(R,G)}$
are locally nilpotent, so we can define
$n_\mu=n_{E_\mu}
\in{\rm Aut}({\frak g}_{{\cal D}}^{\Gamma(R,G)})$ in the
same way as in
(\ref{eqn:nmu}). For each $\beta\in\Pi\setminus\Pi_{{\rm
{max}}}$,
let $E_{\pm
\beta^*}:=n_{\al_\beta}^{-1}n_{\al^*_\beta}E_{\pm \beta}
\in{\frak g}_{{\cal D}}^{\Gamma(R,G)}$,
where we notice that the  $E_{\pm \beta^*}$ are also
locally nilpotent
and that the equalities
$n_{\al_\beta}E_{\pm \beta^*}=n_{\al^*_\beta}E_{\pm
\beta}$
are the same as the ones in ({\rm SR}6-7).
Then we see that the $\Xi$ is surjective since
$\Xi(E_\mu)=E_\mu$ for all $\mu\in{\cal B}$.
To show that the $\Xi$ is injective, it suffices to show
that the elements $h_\sigma$ ($\sigma\in{\cal E}$),
$E_\mu$ ($\mu\in {\cal B}$) of ${\frak g}_{{\cal
D}}^{\Gamma(R,G)}$
satisfy the equalities in ({\rm SR}1-4,$5^\prime$,6-9).
As mentioned above, the equalities in
({\rm SR}6-7)
are satisfied. We only need to check that
the equalities in ({\rm SR}$5^\prime$) are satisfied.
This can be shown as follows. (See Subsec.~3.2 for the
notation $\sim$.)
\par
(1) Assume that $S=\{\al,\beta\}\subset\Pi$ with
$J((\al^*)^\vee,\beta)=-1$ 
and $\Gamma(R,G;S)=\{\al,\al^*,\beta\}$.
Notice that $\al_\beta=\al$. Then we have
\begin{eqnarray*}
\lefteqn{
[E_{\mp\al},E_{\pm\beta^*}]\sim
n_\al^{-1}[E_{\pm\al},n_{\al^*}E_{\pm\beta}]} \\
& &
\sim
n_\al^{-1}[E_{\pm\al},[E_{\pm\al^*},E_{\pm\beta}]]=0\,(\mbox{by
({\rm TSR}10-11)}),
\end{eqnarray*} and
\begin{eqnarray*}
\lefteqn{[E_{\pm\beta},E_{\pm\beta^*}]} \\
& & \sim
n_\al^{-1}[({\rm ad}E_{\pm\al})^{c(\al)}E_{\pm\beta},
[E_{\pm\al^*},E_{\pm\beta}]] \\
& & =
n_\al^{-1}({\rm ad}E_{\pm\al})^{c(\al)-1}(
[[E_{\pm\al},E_{\pm\beta}],
[E_{\pm\al^*},E_{\pm\beta}]]) \,(\mbox{by ({\rm
TSR}10-11)})  \\
& & \sim\left\{\begin{array}{ll}
0 \,(\mbox{by ({\rm TSR}11)})  & \mbox{if $c(\al)=1$,} \\
n_\al^{-1}[E_{\pm\al},
n_\beta [E_{\pm\al},
[E_{\pm\al^*},E_{\pm\beta}]]]=0 \,(\mbox{by ({\rm
TSR}10)}) &
  \mbox{if $c(\al)=2$,}
\end{array}\right. 
\end{eqnarray*}
Since the $E_{\pm\mu}$'s ($\mu\in S\cup S^*$) are locally
nilpotent,
they also satisfy the equalities in ({\rm SR}$5^\prime$)
other than the above ones.
\par
(2) Assume that $S=\{\al,\beta\}\subset\Pi$ with
$J(\beta^\vee,\al)=-1$, $J(\al^\vee,\beta)=-2$, 
$c(\al)=1$ and $\Gamma(R,G;S)=\{\al,\al^*,\beta\}$.
Notice that $\al_\beta=\al$ and $2k(\al)=k(\beta)$. Then,
by ({\rm TSR}8-9),
we have
\begin{equation}\nonumber
[E_{\mp\al},E_{\pm\beta^*}]\sim
n_\al^{-1}[E_{\pm\al},[E_{\pm\al^*},[E_{\pm\al^*},E_{\pm\beta}]]]=0,
\end{equation} and

\begin{eqnarray*}
\lefteqn{[E_{\pm\beta},E_{\pm\beta^*}]
\sim
n_\al^{-1}[({\rm ad}E_{\pm\al})^2E_{\pm\beta},
({\rm ad}E_{\pm\al^*})^2E_{\pm\beta}]} \\
& & = n_\al^{-1}({\rm ad}E_{\pm\al^*})^2
[({\rm ad}E_{\pm\al})^2E_{\pm\beta},E_{\pm\beta}] \\
& & = n_\al^{-1}({\rm ad}E_{\pm\al^*})^2
[[E_{\pm\al},E_{\pm\beta}],[E_{\pm\al},E_{\pm\beta}]]=0.
\end{eqnarray*}
Then we can use the same argument as in (1). \par
(3) Assume that $S=\{\al,\beta,\gamma\}\subset\Pi$ with
$J(\al^\vee,\beta)<0$, $J(\gamma^\vee,\beta)<0$,
$J(\al^\vee,\gamma)=0$, and assume that
$\Gamma(R,G;S)=\{\al,\beta,\beta^*,\gamma\}$
or $\{\al,\al^*,\beta,\beta^*,\gamma\}$.
Then we can use the same argument as in \cite[(2)-(3) of
Proof of
Proposition~4.2]{yamane2}.
\hfill $\BBox$
\section*{Appendix}
As mentioned in the text, especially in
Theorem~\ref{theorem:saitotheorem},
K.~Saito~\cite{saito1} (see also \cite{saitoyoshii})
introduced the notion
of the ERS,
and showed that every RMERS is realized as the $R(k,0)$
for some SEBS
$\EPiako$.
However, there exists a reduced ERS
which is not realized as an RMERS.
As mentioned in Introduction, the authors of
\cite{Allisonetal}
introduced the notion of extended affine root systems
(EARS for short), which is different from the reduced
SEARS's introduced in
\cite{saito1}. It is known (see \cite{Azam}) that there
exists a natural
one-to-one
correspondence between the reduced SEARS's and the EARS's.
Here we also use the terminology and notation in
\cite[Construction~2.32 and Theorem~2.37]{Allisonetal}.
By
Theorems~\ref{theorem:general} and
\ref{theorem:submain}, we see that if an EARS has the
nullity equal to two
and
the rank equal to or more than $2$,  the corresponding
reduced ERS is
realized as
$R(k,g)$ for some QEBS ${\cal D}=\EPiakg$ with
$g(\al)=\emptyset$ or
$2\BBbb{Z}+1$
($\al\in\Pi$).
Let ${\cal D}$ be a QEBS with $l\geq 2$ such that the name
of $A(=A_\Pi)$
is ${\bf D}_{l+1}^{(2)}$.
Notice that there exist $\varepsilon_i\in{\cal E}$
($1\leq i\leq l$) such that
$J(\varepsilon_i,\varepsilon_j)=\delta_{ij}$,
$\al_0=\delta-\varepsilon_1$,
$\al_i=\varepsilon_i-\varepsilon_{i+1}$
($2\leq i\leq l-1$) and $\al_l=\varepsilon_l$.
Then the corresponding EARS is $R(X,S,L,E)$
is such that
\begin{equation} \nonumber
X=\left\{\begin{array}{ll}
B_l & \mbox{if $g(\al_0)=g(\al_l)=\emptyset$} \\
BC_l & \mbox{otherwise,}
\end{array}\right.
\end{equation}
and $S=((2\BBbb{Z}+1)\delta+\BBbb{Z}k(\al_0)a)
\cup(2\BBbb{Z}\delta+\BBbb{Z}k(\al_l)a))$,
$L=2\BBbb{Z}\delta+\BBbb{Z}k(\al_i)a$ ($2\leq i\leq l-1$)
and $E=((4\BBbb{Z}+2)\delta+g(\al_0)a)
\cup(4\BBbb{Z}\delta+g(\al_l)a))$.
Here if $g(\al_0)=\emptyset$, then
$(4\BBbb{Z}+2)\delta+g(\al_0)a)=\emptyset$;
if $g(\al_1)=\emptyset$, then
$(4\BBbb{Z}+2)\delta+g(\al_1)a)=\emptyset$.
(Strictly speaking, in \cite{Allisonetal},
$R(X,S,L,\emptyset)$ is denoted as $R(X,S,L)$).

\end{document}